# Some integrals involving the Stieltjes constants: Part II

Donal F. Connon

dconnon@btopenworld.com

11 April 2011

**Abstract**

Some new integrals involving the Stieltjes constants are developed in this paper.

## 1. Introduction

The Stieltjes constants $\gamma_p(x)$ are the coefficients of the Laurent expansion of the Hurwitz zeta function $\varsigma(s,x)$ about $s=1$

$$(1.1) \qquad \varsigma(s,x) = \sum_{n=0}^{\infty} \frac{1}{(n+x)^s} = \frac{1}{s-1} + \sum_{p=0}^{\infty} \frac{(-1)^p}{p!} \gamma_p(x)(s-1)^p$$

where $\gamma_p(x)$ are known as the generalised Stieltjes constants and we have [41]

$$(1.2) \qquad \gamma_0(x) = -\psi(x)$$

where $\psi(x)$ is the digamma function.

With $x=1$ equation (1.1) reduces to the Riemann zeta function

$$\varsigma(s) = \frac{1}{s-1} + \sum_{p=0}^{\infty} \frac{(-1)^p}{p!} \gamma_p(s-1)^p$$

As previously noted in [20], using (1.1) it is easily seen that the difference of two Stieltjes constants may be represented by

(1.3) $$\gamma_p(x) - \gamma_p(y) = (-1)^p \lim_{s \to 1} \frac{\partial^p}{\partial s^p}[\varsigma(s,x) - \varsigma(s,y)]$$

## 2. A useful formula for the Stieltjes constants

We recall Hasse's formula [31] for the Hurwitz zeta function which is valid for all $s \in \mathbb{C}$ except $s = 1$ (in this form, it is valid in the limit as $s \to 1$)

(2.1) $$(s-1)\varsigma(s,x) = \sum_{j=0}^{\infty} \frac{1}{j+1} \sum_{k=0}^{j} \binom{j}{k} \frac{(-1)^k}{(x+k)^{s-1}}$$

and differentiation with respect to $x$ gives us

$$\frac{\partial}{\partial x}\varsigma(s,x) = -\sum_{j=0}^{\infty} \frac{1}{j+1} \sum_{k=0}^{j} \binom{j}{k} \frac{(-1)^k}{(x+k)^s}$$

and we then have

(2.2) $$\frac{\partial^n}{\partial s^n} \frac{\partial}{\partial x} \varsigma(s,x) = (-1)^{n+1} \sum_{j=0}^{\infty} \frac{1}{j+1} \sum_{k=0}^{j} \binom{j}{k} \frac{(-1)^k \log^n(x+k)}{(x+k)^s}$$

We note that the partial derivatives commute in the region where $\varsigma(s,x)$ is analytic and hence we have

$$\frac{\partial^n}{\partial s^n} \frac{\partial}{\partial x} \varsigma(s,x) = \frac{\partial}{\partial x} \frac{\partial^n}{\partial s^n} \varsigma(s,x)$$

Evaluation of (2.2) at $s = 0$ results in

(2.3) $$\frac{\partial}{\partial x} \varsigma^{(n)}(0,x) = (-1)^{n+1} \sum_{j=0}^{\infty} \frac{1}{j+1} \sum_{k=0}^{j} \binom{j}{k} (-1)^k \log^n(x+k)$$

We may write (1.1) as

(2.4) $$\varsigma(s+1,x) = \frac{1}{s} + \sum_{n=0}^{\infty} \frac{(-1)^n}{n!} \gamma_n(x) s^n$$

and we have the Maclaurin expansion

(2.5) $$\varsigma(s,x) = \frac{1}{s-1} + \sum_{n=0}^{\infty} \frac{(-1)^{n+1} R_n(x)}{n!} s^n$$

where

$$R_n(x) = (-1)^{n+1} \left. \frac{\partial^n}{\partial s^n} \varsigma(s,x) \right|_{s=0}$$

$R(x) = R_2(x)$ is referred to as the Deninger $R$-function, after Deninger [28] who introduced it in 1984.

Differentiating (2.5) with respect to $x$ gives us

(2.6) $$\frac{\partial}{\partial x} \varsigma(s,x) = \sum_{n=0}^{\infty} \frac{(-1)^{n+1} R_n'(x)}{n!} s^n$$

We note that

$$\frac{\partial}{\partial x} \varsigma(s,x) = -s\varsigma(s+1,x)$$

and, using (2.4), this is equal to

(2.7) $$= -1 + \sum_{n=0}^{\infty} \frac{(-1)^{n+1}}{n!} \gamma_n(x) s^{n+1}$$

and comparing the coefficients of (2.6) and (2.7) Chakraborty, Kanemitsu and Kuzumaki [12] deduced the important identity

(2.8) $$R_n'(x) = (-1)^{n+1} \left. \frac{\partial^n}{\partial s^n} \frac{\partial}{\partial x} \varsigma(s,x) \right|_{s=0} = -n\gamma_{n-1}(x)$$

This may be compared with the more familiar formula for the Stieltjes constants where the limit is evaluated at $s=1$

(2.9) $$(-1)^{n+1} \left. \frac{\partial^{n-1}}{\partial s^{n-1}} \left[ \varsigma(s,x) - \frac{1}{s-1} \right] \right|_{s=1} = -\gamma_{n-1}(x)$$

We see from (2.3) that (2.8) is equal to

$$= \sum_{j=0}^{\infty} \frac{1}{j+1} \sum_{k=0}^{j} \binom{j}{k} (-1)^k \log^n(x+k)$$

and hence we easily deduce that

$$(2.10) \qquad \gamma_n(x) = -\frac{1}{n+1}\sum_{j=0}^{\infty}\frac{1}{j+1}\sum_{k=0}^{j}\binom{j}{k}(-1)^k \log^{n+1}(x+k)$$

which was previously obtained in 2007 in [20].

In [20, Eq. (4.3.231)] we also noted that $R'_n(x) = -n\gamma_{n-1}(x)$ in the equivalent form for $x > 0$

$$\int_1^x \gamma_n(t)\,dt = \frac{(-1)^{n+1}}{n+1}\left[\varsigma^{(n+1)}(0,x) - \varsigma^{(n+1)}(0)\right]$$

but the usefulness of this simple identity was not then fully appreciated by the author.

The formula (2.8) features throughout the rest of this paper where it is used to simplify the derivation of some known identities and also to produce some new ones.

**Remark (i)**

The equivalent formula to (2.5) as reported by Chakraborty et al. [12] did not include the term $\frac{1}{s-1}$; it seems to me that it should be so included if only to concur with the analysis previously carried out by Sitaramachandrarao [37] in 1986 where he considered the Maclaurin series for the Riemann zeta function

$$(2.11) \qquad \varsigma(s) + \frac{1}{1-s} = \sum_{n=0}^{\infty}\frac{(-1)^n \delta_n}{n!}s^n$$

where

$$(2.12) \qquad \delta_n = \lim_{m\to\infty}\left[\sum_{k=1}^{m}\log^n k - \int_1^m \log^n x\,dx - \frac{1}{2}\log^n m\right]$$

$$= (-1)^n\left[\varsigma^{(n)}(0) + n!\right]$$

This additional term $\frac{1}{s-1}$ of course vanishes when equation (2.5) is differentiated with respect to $x$ and thus (2.8) continues to remain valid.

**Remark (ii)**

We note from (2.8) that with $n=1$

$$R_1'(x) = \frac{\partial}{\partial s}\frac{\partial}{\partial x}\varsigma(s,x)\Big|_{s=0} = -\gamma_0(x) = \psi(x)$$

or equivalently

$$\frac{\partial}{\partial x}\varsigma'(0,x) = \psi(x)$$

and integration then results in

(2.13) $\qquad \varsigma'(0,t) - \varsigma'(0) = \log \Gamma(t)$

We have Legendre's duplication formula for the gamma function [38, p.7]

$$\log \Gamma(2t) = \log \Gamma(t) + \log \Gamma\left(t+\frac{1}{2}\right) + (2t-1)\log 2 - \frac{1}{2}\log \pi$$

and substituting (2.13) gives us

$$\varsigma'(0,2t) = \varsigma'(0,t) + \varsigma'\left(0,t+\frac{1}{2}\right) - \varsigma'(0) + (2t-1)\log 2 - \frac{1}{2}\log \pi$$

and with $t = 1/2$ we obtain

$$\varsigma'\left(0,\frac{1}{2}\right) = \varsigma'(0) + \frac{1}{2}\log \pi$$

We recall the identity [30]

$$\varsigma\left(s,\frac{1}{2}\right) = [2^s - 1]\varsigma(s)$$

and differentiation results in

$$\varsigma'\left(s,\frac{1}{2}\right) = [2^s - 1]\varsigma'(s) + \varsigma(s)2^s \log 2$$

so that

$$\varsigma'\left(0,\frac{1}{2}\right) = \varsigma(0)\log 2 = -\frac{1}{2}\log 2$$

This then gives us the well-known result

(2.14) $\qquad \varsigma'(0) = -\frac{1}{2}\log(2\pi)$

and hence we have obtained Lerch's identity [7] in a very direct manner without the need to resort to the functional equation for the Riemann zeta function

(2.15) $\qquad \varsigma'(0,t) = \log\Gamma(t) - \frac{1}{2}\log(2\pi)$

**Remark (iii)**

It should be noted that the Stieltjes constant $\gamma_1$ referred to in Eq. (2.15) of Deninger's paper [28] should be increased by a factor of 2 (this difference arises because Deninger [28, p.174] employed a different definition in the Laurent expansion of the Hurwitz zeta function; this value was also inconsistently employed in Eq. (1.26) of the paper by Chakraborty et al. [12]).

**3. An application of the Abel-Plana summation formula**

Adamchik [2] has recently reported that the Hermite integral for the Hurwitz zeta function may be derived from the Abel-Plana summation formula [38, p.90]

(3.1) $\qquad \sum_{k=0}^{\infty} f(k) = \frac{1}{2}f(0) + \int_0^{\infty} f(x)\,dx + i\int_0^{\infty} \frac{f(ix) - f(-ix)}{e^{2\pi x} - 1}\,dx$

which applies to functions which are analytic in the right-hand plane and satisfy the convergence condition $\lim_{y \to \infty} e^{-2\pi y}|f(x+iy)| = 0$ uniformly on any finite interval of $x$. Derivations of the Abel-Plana summation formula are given in [39, p.145] and [40, p.108].

The Hermite integral for the Hurwitz zeta function may be derived as follows. Letting $f(k) = (k+u)^{-s}$ we obtain

(3.2) $\qquad \varsigma(s,u) = \sum_{k=0}^{\infty} \frac{1}{(k+u)^s} = \frac{u^{-s}}{2} + \frac{u^{1-s}}{s-1} + i\int_0^{\infty} \frac{(u+ix)^{-s} - (u-ix)^{-s}}{e^{2\pi x} - 1}\,dx$

Then, noting that

$$(u+ix)^{-s} - (u-ix)^{-s} = (re^{i\theta})^{-s} - (re^{-i\theta})^{-s}$$

$$= r^{-s}[e^{-is\theta} - e^{is\theta}]$$

(3.2.1) $$(u+ix)^{-s} - (u-ix)^{-s} = \frac{2}{i(u^2+x^2)^{s/2}} \sin(s \tan^{-1}(x/u))$$

we may write (3.2) as Hermite's integral for the Hurwitz zeta function $\varsigma(s,u)$

(3.3) $$\varsigma(s,u) = \frac{u^{-s}}{2} + \frac{u^{1-s}}{s-1} + 2\int_0^\infty \frac{\sin(s \tan^{-1}(x/u))}{(u^2+x^2)^{s/2}(e^{2\pi x}-1)} dx$$

Differentiating (3.2) with respect to $u$

(3.4) $$\frac{\partial}{\partial u}\varsigma(s,u) = -\frac{su^{-s-1}}{2} - u^{-s} - is\int_0^\infty \frac{(u+ix)^{-s-1} - (u-ix)^{-s-1}}{e^{2\pi x}-1} dx$$

We note that if $f(s) = sg(s)$ then the Leibniz differentiation formula results in

$$f^{(n+1)}(s) = s\, g^{(n+1)}(s) + (n+1)g^{(n)}(s)$$

so that

(3.5) $$f^{(n+1)}(0) = (n+1)g^{(n)}(0)$$

and hence we obtain

$$\frac{\partial^{n+1}}{\partial s^{n+1}} \frac{\partial}{\partial u}\varsigma(s,u)\bigg|_{s=0} = -\frac{n+1}{2u}(-1)^n \log^n u - (-1)^{n+1}\log^{n+1} u$$

$$-i(-1)^n(n+1)\int_0^\infty \frac{(u-ix)\log^n(u+ix) - (u+ix)\log^n(u-ix)}{(u^2+x^2)(e^{2\pi x}-1)} dx$$

and comparing this with (2.8) we readily obtain the integral formula originally obtained by Coffey [17] in 2007 for the Stieltjes constants

(3.6) $$\gamma_n(u) = \frac{1}{2u}\log^n u - \frac{1}{n+1}\log^{n+1} u + i\int_0^\infty \frac{(u-ix)\log^n(u+ix) - (u+ix)\log^n(u-ix)}{(u^2+x^2)(e^{2\pi x}-1)} dx$$

This derivation is slightly more direct than the one originally provided by Coffey [17] (and is also simpler than my previous proof in [22]).

□

Chen [13] has recently shown that for $u > 0$ and $s > 0$

(3.7) $$\frac{2}{\Gamma(s)}\int_0^\infty e^{-uy^2} y^{2s-1} \sin(xy^2)\,dy = \frac{\sin[s\tan^{-1}(x/u)]}{(u^2+x^2)^{s/2}}$$

and hence we have using Hermite's integral (3.3)

$$\varsigma(s,u) = \frac{u^{-s}}{2} + \frac{u^{1-s}}{s-1} + \frac{4}{\Gamma(s)}\int_0^\infty \frac{1}{e^{2\pi x}-1}\,dx \int_0^\infty e^{-uy^2} y^{2s-1} \sin(xy^2)\,dy$$

$$= \frac{u^{-s}}{2} + \frac{u^{1-s}}{s-1} + \frac{4}{\Gamma(s)}\int_0^\infty e^{-uy^2} y^{2s-1}\,dy \int_0^\infty \frac{\sin(xy^2)}{e^{2\pi x}-1}\,dx$$

Using Legendre's relation [39, p.122]

(3.7.1) $$2\int_0^\infty \frac{\sin(xt)}{e^{2\pi x}-1}\,dx = \frac{1}{e^t-1} - \frac{1}{t} + \frac{1}{2} = \frac{1}{2}\coth\frac{t}{2} - \frac{1}{t}$$

this results in

$$\varsigma(s,u) = \frac{u^{-s}}{2} + \frac{u^{1-s}}{s-1} + \frac{2}{\Gamma(s)}\int_0^\infty e^{-uy^2} y^{2s-1}\left[\frac{1}{e^{y^2}-1} - \frac{1}{y^2} + \frac{1}{2}\right]dy$$

With the substitution $v = y^2$ this becomes

(3.8) $$\varsigma(s,u) = \frac{u^{-s}}{2} + \frac{u^{1-s}}{s-1} + \frac{1}{\Gamma(s)}\int_0^\infty e^{-uv} v^{s-1}\left[\frac{1}{e^v-1} - \frac{1}{v} + \frac{1}{2}\right]dv$$

which is reported in [38, p.92] as being valid for $\operatorname{Re}(s) > -1$.

With $u = 1$ we have

(3.9) $$\varsigma(s) = \frac{1}{2} + \frac{1}{s-1} + \frac{1}{\Gamma(s)}\int_0^\infty e^{-v} v^{s-1}\left[\frac{1}{e^v-1} - \frac{1}{v} + \frac{1}{2}\right]dv$$

which is also reported in [38, p.100] as being valid for $\operatorname{Re}(s) > -1$.

$\square$

We may write (3.7) as

$$\frac{2}{\Gamma(s)}\int_0^\infty e^{-uy^2} y^{2s-1} \sin(xy^2)\,dy = \frac{2s}{\Gamma(1+s)}\int_0^\infty e^{-uy^2} y^{2s-1} \sin(xy^2)\,dy$$

and we have the limit

$$\lim_{s \to 0} \frac{2}{\Gamma(s)} \int_0^\infty e^{-uy^2} y^{2s-1} \sin(xy^2) \, dy = \lim_{s \to 0} \frac{2s}{\Gamma(1+s)} \int_0^\infty e^{-uy^2} y^{2s-1} \sin(xy^2) \, dy$$

$$= \lim_{s \to 0} \frac{2s}{\Gamma(1+s)} \int_0^\infty \frac{e^{-uy} \sin(xy)}{y} \, dy$$

We have the well-known integral (rigorous derivations of which are contained in [3, p.285] and [6, p.272])

(3.10) $$\tan^{-1}(x/u) = \int_0^\infty \frac{e^{-uy} \sin(xy)}{y} \, dy$$

and hence we see that

$$\lim_{s \to 0} \frac{2}{\Gamma(s)} \int_0^\infty e^{-uy^2} y^{2s-1} \sin(xy^2) \, dy = 0$$

showing that Chen's result is valid for $s \geq 0$.

□

For completeness we present another proof of (3.7). Using the definition of the gamma function we have

$$\int_0^\infty e^{-zy} y^{s-1} \, dy = \frac{\Gamma(s)}{z^s}$$

and with $z = u \pm ix$ we obtain

$$\int_0^\infty e^{-uy} y^{s-1} [\cos(xy) - i\sin(xy)] \, dy = (u+ix)^{-s} \Gamma(s)$$

$$\int_0^\infty e^{-uy} y^{s-1} [\cos(xy) + i\sin(xy)] \, dy = (u-ix)^{-s} \Gamma(s)$$

This gives us

$$i \int_0^\infty e^{-uy} y^{s-1} \sin(xy) \, dy = -\left[ (u+ix)^{-s} - (u-ix)^{-s} \right] \Gamma(s)$$

$$\int_0^\infty e^{-uy} y^{s-1} \cos(xy)\, dy = \left[(u+ix)^{-s} + (u-ix)^{-s}\right]\Gamma(s)$$

Then, noting (3.2.1) we obtain

$$\frac{1}{\Gamma(s)}\int_0^\infty e^{-uy} y^{s-1} \sin(xy)\, dy = \frac{\sin[s\tan^{-1}(x/u)]}{(u^2+x^2)^{s/2}}$$

and

$$\frac{1}{\Gamma(s)}\int_0^\infty e^{-uy} y^{s-1} \cos(xy)\, dy = \frac{\cos[s\tan^{-1}(x/u)]}{(u^2+x^2)^{s/2}}$$

$\square$

We have from (3.8)

$$(3.11) \quad \Gamma(s)\left[\varsigma(s,u) - \frac{u^{-s}}{2} - \frac{u^{1-s}}{s-1}\right] = \int_0^\infty e^{-uv} v^{s-1}\left[\frac{1}{e^v-1} - \frac{1}{v} + \frac{1}{2}\right] dv$$

and we consider the limit as $s \to 0$. With simple algebra we may write

$$\varsigma(s,u) - \frac{u^{-s}}{2} - \frac{u^{1-s}}{s-1} = \varsigma(s,u) - \varsigma(0,u) + \left(\frac{1}{2} - u\right) - \frac{u^{-s}}{2} - \frac{u^{1-s}}{s-1}$$

$$= \varsigma(s,u) - \varsigma(0,u) + \left(\frac{1}{2} - u\right)(1 - u^{-s}) - u^{1-s}\left(1 + \frac{1}{s-1}\right)$$

and we have

$$\Gamma(s)\left[\varsigma(s,u) - \frac{u^{-s}}{2} - \frac{u^{1-s}}{s-1}\right] = \Gamma(1+s)\left[\frac{\varsigma(s,u) - \varsigma(0,u)}{s} + \left(\frac{1}{2} - u\right)\frac{(1 - u^{-s})}{s} - \frac{u^{1-s}}{s-1}\right]$$

Taking the limit as $s \to 0$ we see that

$$\lim_{s \to 0} \Gamma(s)\left[\varsigma(s,u) - \frac{u^{-s}}{2} - \frac{u^{1-s}}{s-1}\right] = \varsigma'(0,u) + \left(\frac{1}{2} - u\right)\log u + u$$

and hence we obtain

$$(3.12) \quad \varsigma'(0,u) + \left(\frac{1}{2} - u\right)\log u + u = \int_0^\infty \frac{e^{-uv}}{v}\left[\frac{1}{e^v-1} - \frac{1}{v} + \frac{1}{2}\right] dv$$

Using Lerch's identity (2.15), it may be seen that this is equivalent to Binet's first formula for $\log \Gamma(u)$ [39, p.249]

$$(3.13) \quad \log \Gamma(u) = \left(u - \frac{1}{2}\right)\log u - u + \frac{1}{2}\log(2\pi) + \int_0^\infty \frac{e^{-uy}}{y}\left[\frac{1}{e^y - 1} - \frac{1}{y} + \frac{1}{2}\right]dy$$

We note Alexeiewsky's theorem [38, p.32]

$$(3.14) \quad \int_0^t \log \Gamma(u)du = \frac{1}{2}t(1-t) + \frac{1}{2}t\log(2\pi) - \log G(1+t) + t\log \Gamma(t)$$

where $G(x)$ is the Barnes double gamma function $\Gamma_2(x) = 1/G(x)$ defined, inter alia, by the Weierstrass canonical product [38, p.25]

$$(3.15) \quad G(1+x) = (2\pi)^{x/2} \exp\left[-\frac{1}{2}(\gamma x^2 + x^2 + x)\right]\prod_{k=1}^\infty \left\{\left(1 + \frac{x}{k}\right)^k \exp\left(\frac{x^2}{2k} - x\right)\right\}$$

and we note that $G(1) = G(2) = 1$.

Hence integrating (3.13) results in

$$(3.16) \quad \log G(1+t) = t\log \Gamma(t) + \frac{1}{4}t^2 - \frac{1}{2}t(t-1)\log t - \int_0^\infty \frac{1-e^{-tv}}{v^2}\left[\frac{1}{e^v - 1} - \frac{1}{v} + \frac{1}{2}\right]dv$$

and with $t = 1$ we immediately obtain

$$\int_0^\infty \frac{1-e^{-v}}{v^2}\left[\frac{1}{e^v - 1} - \frac{1}{v} + \frac{1}{2}\right]dv = \frac{1}{4}$$

Subtracting this from (3.16) results in

$$(3.17) \quad \log G(1+t) = t\log \Gamma(t) + \frac{1}{4}(t^2 - 1) - \frac{1}{2}t(t-1)\log t + \int_0^\infty \frac{e^{-tv} - e^{-v}}{v^2}\left[\frac{1}{e^v - 1} - \frac{1}{v} + \frac{1}{2}\right]dv$$

It may be possible employ Pringsheim's artifice with this integral (as was employed in [39, p.249]).

□

We obtain by differentiating (3.2)

$$\varsigma'(0,u) = \left(u - \frac{1}{2}\right)\log u - u + 2\int_0^\infty \frac{\tan^{-1}(x/u)}{e^{2\pi x} - 1} dx$$

and differentiating this with respect to $u$ gives us

$$\frac{\partial}{\partial u}\varsigma'(0,u) = \log u - \frac{1}{2u} - 2\int_0^\infty \frac{x}{(u^2 + x^2)(e^{2\pi x} - 1)} dx$$

Since from (2.8) $\frac{\partial}{\partial u}\varsigma'(0,u) = -\gamma_0(u)$ we obtain

$$-\gamma_0(u) = \log u - \frac{1}{2u} - 2\int_0^\infty \frac{x}{(u^2 + x^2)(e^{2\pi x} - 1)} dx$$

and since $\gamma_0(u) = -\psi(u)$ this is equivalent to

$$\psi(u) = \log u - \frac{1}{2u} - 2\int_0^\infty \frac{x}{(u^2 + x^2)(e^{2\pi x} - 1)} dx$$

as reported in [39, p.251].

Similarly, we have

$$\varsigma''(0,u) = \left(\frac{1}{2} - u\right)\log^2 u + 2u \log u - 2u - 2\int_0^\infty \frac{\log(u^2 + x^2)\tan^{-1}(x/u)}{e^{2\pi x} - 1} dx$$

so that

$$\frac{\partial}{\partial u}\varsigma''(0,u) = \frac{1}{u}\log u - \log^2 u$$

$$-2\pi\int_0^\infty \frac{x}{(u^2 + x^2)(e^{2\pi x} - 1)} dx + 2\int_0^\infty \frac{x\log(u^2 + x^2)}{(u^2 + x^2)(e^{2\pi x} - 1)} dx + 4u\int_0^\infty \frac{\tan^{-1}(u/x)}{(u^2 + x^2)(e^{2\pi x} - 1)} dx$$

Then using $\left.\frac{\partial^2}{\partial s^2}\frac{\partial}{\partial u}\varsigma(s,u)\right|_{s=0} = 2\gamma_1(u)$ we obtain

$$\text{(3.18)} \quad \gamma_1(u) = \frac{1}{2u}\log u - \frac{1}{2}\log^2 u + \int_0^\infty \frac{x\log(u^2+x^2)}{(u^2+x^2)(e^{2\pi x}-1)}dx - 2u\int_0^\infty \frac{\tan^{-1}(x/u)}{(u^2+x^2)(e^{2\pi x}-1)}dx$$

It may be noted that Shail [36, p.799] made reference to these "seemingly intractable" integrals in 2000.

□

We have the well-known Hurwitz's formula for the Fourier expansion of the Hurwitz zeta function $\varsigma(s,t)$ as reported in Titchmarsh's treatise

$$\varsigma(s,t) = 2\Gamma(1-s)\left[\sin\left(\frac{\pi s}{2}\right)\sum_{n=1}^\infty \frac{\cos 2n\pi t}{(2\pi n)^{1-s}} + \cos\left(\frac{\pi s}{2}\right)\sum_{n=1}^\infty \frac{\sin 2n\pi t}{(2\pi n)^{1-s}}\right]$$

where $\operatorname{Re}(s) < 0$ and $0 < t \leq 1$. In 2000, Boudjelkha showed that this formula also applies in the region $\operatorname{Re}(s) < 1$. It may be noted that when $t = 1$ this reduces to Riemann's functional equation for $\varsigma(s)$.

Unfortunately, it appears that (2.8) cannot be applied to Hurwitz's formula for the Fourier expansion of the Hurwitz zeta function because this would result in divergent series. This however may give an indication of Ramanujan's erroneous thinking in this area as exemplified in Berndt's book, Ramanujan's Notebooks, Part I [8, p.200].

□

It is instructive to consider the original derivation of the formula for the Stieltjes constants which was derived by Stieltjes himself. This is recorded in a letter dated June 1885 from Stieltjes to Hermite [32] and the proof proceeds as follows:

From (3.9) we have

$$\varsigma(1+s) - \frac{1}{s} = \frac{1}{2} + \frac{1}{\Gamma(1+s)}\int_0^\infty e^{-x}\left[\frac{1}{e^x-1} - \frac{1}{x} + \frac{1}{2}\right]x^s dx$$

$$= \frac{1}{2} + \frac{1}{\Gamma(1+s)}\int_0^\infty e^{-x}\left[\frac{e^{-x}}{1-e^{-x}} - \frac{1}{x} + \frac{1}{2}\right]x^s dx$$

$$= \frac{1}{2} + \frac{1}{\Gamma(1+s)}\int_0^\infty e^{-x}\left[\frac{1}{1-e^{-x}} - \frac{1}{x} - \frac{1}{2}\right]x^s dx$$

so that

$$\text{(3.19)} \quad \varsigma(1+s) - \frac{1}{s} = \frac{1}{\Gamma(1+s)}\int_0^\infty e^{-x}\left[\frac{1}{1-e^{-x}} - \frac{1}{x}\right]x^s dx$$

Using the Maclaurin expansion

$$x^s = \sum_{j=0}^{\infty} \frac{\log^j x}{j!} s^j$$

we see that

$$\varsigma(1+s) - \frac{1}{s} = \frac{1}{\Gamma(1+s)} \sum_{j=0}^{\infty} \frac{a_j}{j!} s^j$$

where $a_j = \int_0^{\infty} e^{-x} \left[ \frac{1}{1-e^{-x}} - \frac{1}{x} \right] \log^j x \, dx$.

Let $F_n(t) = \int_0^{\infty} e^{-tx} \log^n x \, dx$. We have $\Gamma^{(n)}(1) = \int_0^{\infty} e^{-x} \log^n x \, dx$ and with $x \to tx$ this becomes

$$\Gamma^{(n)}(1) = t \int_0^{\infty} e^{-tx} (\log x + \log t)^n \, dx$$

We write $F_n(t)$ as

$$F_n(t) = \int_0^{\infty} e^{-tx} (\log x + \log t - \log t)^n \, dx$$

$$= \frac{1}{t} \sum_{j=0}^{n} (-1)^j \binom{n}{j} \Gamma^{(n-j)}(1) \log^j t$$

and we then have

(3.20) $$\sum_{j=0}^{n} (-1)^j \binom{n}{j} \Gamma^{(n-j)}(1) \frac{\log^j t}{t} = \int_0^{\infty} e^{-tx} \log^n x \, dx$$

We have the summation

$$\sum_{j=0}^{n} (-1)^j \binom{n}{j} \Gamma^{(n-j)}(1) \sum_{t=1}^{r} \frac{\log^j t}{t} = \sum_{t=1}^{r} \int_0^{\infty} e^{-tx} \log^n x \, dx$$

$$= \int_0^{\infty} \left[ \frac{e^{-x} - e^{-(r+1)x}}{1 - e^{-x}} \right] \log^n x \, dx$$

and the right-hand side may be written as $a_n + f_n(r) - g_n(r)$ where

$$(3.21) \qquad g_n(r) = \int_0^\infty \left[\frac{1}{1-e^{-x}} - \frac{1}{x}\right] e^{-(r+1)x} \log^n x \, dx$$

$$(3.22) \qquad f_n(r) = \int_0^\infty \left[\frac{e^{-x}}{x} - \frac{e^{-(r+1)x}}{x}\right] \log^n x \, dx$$

where we note that $\lim_{r \to \infty} g_n(r) = 0$. Now treating $r$ as a continuous variable we obtain

$$f_n'(r) = F_n(r+1) = \sum_{j=0}^n (-1)^j \binom{n}{j} \Gamma^{(n-j)}(1) \frac{\log^j(r+1)}{r+1}$$

and since $f_n(0) = 0$ we have by integration

$$(3.23) \qquad f_n(r) = \sum_{j=0}^n (-1)^j \binom{n}{j} \Gamma^{(n-j)}(1) \frac{\log^{j+1}(r+1)}{j+1}$$

so

$$a_n - \int_0^\infty \left[\frac{1}{1-e^{-x}} - \frac{1}{x}\right] e^{-(r+1)x} \log^n x \, dx = \sum_{j=0}^n (-1)^j \binom{n}{j} \Gamma^{(n-j)}(1) \left(\sum_{m=1}^r \frac{\log^j m}{m} - \frac{\log^{j+1}(r+1)}{j+1}\right)$$

$$= \sum_{j=0}^n (-1)^j \binom{n}{j} \Gamma^{(n-j)}(1) \left(\sum_{m=1}^r \frac{\log^j m}{m} - \frac{\log^{j+1} r + \left[\log^{j+1}(r+1) - \log^{j+1} r\right]}{j+1}\right)$$

Therefore letting $r \to \infty$ we obtain

$$a_n = \sum_{j=0}^n (-1)^j \binom{n}{j} \Gamma^{(n-j)}(1) \lim_{r \to \infty} \left(\sum_{m=1}^r \frac{\log^j m}{m} - \frac{\log^{j+1} r}{j+1}\right)$$

since $\lim_{r \to \infty} \left[\log^{j+1}(r+1) - \log^{j+1} r\right] = 0$. This may be written as

$$(3.24) \qquad a_n = \sum_{j=0}^n \binom{n}{j} (-1)^j C_j \, \Gamma^{(n-j)}(1)$$

where

$$C_j = \lim_{r \to \infty} \left(\sum_{m=1}^r \frac{\log^j m}{m} - \frac{\log^{j+1} r}{j+1}\right)$$

Eq. (3.24) reminds us of the Cauchy product of two series

$$\sum_{j=0}^{\infty}\frac{a_j}{j!}s^j = \left(\sum_{m=0}^{\infty}(-1)^m\frac{C_m}{m!}s^j\right)\left(\sum_{k=0}^{\infty}\frac{\Gamma^{(k)}(1)}{k!}s^j\right)$$

where $\Gamma(1+s) = \sum_{k=0}^{\infty}\frac{\Gamma^{(k)}(1)}{k!}s^j$.

We have from (3.19)

(3.25) $$\left[\varsigma(1+s) - \frac{1}{s}\right]\Gamma(1+s) = \int_0^{\infty}e^{-x}\left[\frac{1}{1-e^{-x}} - \frac{1}{x}\right]x^s\,dx$$

and we note from (2.4) that

$$\varsigma(s+1) - \frac{1}{s} = \sum_{m=0}^{\infty}\frac{(-1)^m}{m!}\gamma_m s^m$$

and we conclude that

(3.26) $$\gamma_j = C_j = \lim_{r\to\infty}\left(\sum_{m=1}^{r}\frac{\log^j m}{m} - \frac{\log^{j+1} r}{j+1}\right)$$

A more direct derivation of (3.26) is given, inter alia, by Bohman and Fröberg [10].

## 4. A family of integral representations of the Stieltjes constants

With a view to utilising (2.8), we first of all differentiate (3.8) with respect to $u$ to give us

$$\frac{\partial}{\partial u}\varsigma(s,u) = -\frac{su^{-s-1}}{2} - u^{-s} - \frac{1}{\Gamma(s)}\int_0^{\infty}e^{-uv}v^s\left[\frac{1}{e^v-1} - \frac{1}{v} + \frac{1}{2}\right]dv$$

and differentiating this $n$ times with respect to $s$ (with the assistance of (3.5)) gives us

(4.1)
$$\left.\frac{\partial^n}{\partial s^n}\frac{\partial}{\partial u}\varsigma(s,u)\right|_{s=0} = -\frac{n}{2u}(-1)^{n-1}\log^{n-1}u - (-1)^n\log^n u - \int_0^{\infty}F^{(n)}(0,v)e^{-uv}\left[\frac{1}{e^v-1} - \frac{1}{v} + \frac{1}{2}\right]dv$$

where, for convenience, we have designated $F(s,v)$ as

(4.2) $$F(s,v) = \frac{v^s}{\Gamma(s)} = \frac{sv^s}{\Gamma(1+s)}$$

and we note that $F(0,v) = 0$. We employ the notation $\dfrac{\partial^n}{\partial s^n} F(s,v) = F^{(n)}(s,v)$ and

$$\left. \dfrac{\partial^n}{\partial s^n} F(s,v) \right|_{s=0} = F^{(n)}(0,v).$$

Then using (2.8)

$$\left. \dfrac{\partial^n}{\partial s^n} \dfrac{\partial}{\partial u} \varsigma(s,u) \right|_{s=0} = (-1)^n n \gamma_{n-1}(u)$$

we have thus determined that for $n \geq 1$

(4.3) $\quad \gamma_{n-1}(u) = \dfrac{1}{2u} \log^{n-1} u - \dfrac{1}{n} \log^n u - \dfrac{(-1)^n}{n} \int_0^\infty F^{(n)}(0,v) e^{-uv} \left[ \dfrac{1}{e^v - 1} - \dfrac{1}{v} + \dfrac{1}{2} \right] dv$

and, in particular, for $n \geq 2$ we have a family of integral representations of the Stieltjes constants

(4.4) $\quad \gamma_{n-1} = \gamma_{n-1}(1) = \dfrac{(-1)^{n+1}}{n} \int_0^\infty F^{(n)}(0,v) e^{-uv} \left[ \dfrac{1}{e^v - 1} - \dfrac{1}{v} + \dfrac{1}{2} \right] dv$

and for $n = 1$ we obtain

$$\gamma = \dfrac{1}{2} + \int_0^\infty F^{(1)}(0,v) e^{-v} \left[ \dfrac{1}{e^v - 1} - \dfrac{1}{v} + \dfrac{1}{2} \right] dv$$

We shall initially consider the first two derivatives of $F(s,v)$. Differentiating (4.2) with respect to $s$ results in

$$F'(s,v) = \dfrac{sv^s}{\Gamma(1+s)} \left[ \log v + \dfrac{1}{s} - \psi(1+s) \right]$$

$$= \dfrac{sv^s}{\Gamma(1+s)} [\log v - \psi(s)]$$

and so we have

(4.5) $\quad F'(s,v) = F(s,v)[\log v - \psi(s)]$

Since $s\psi(s) = s\psi(1+s) - 1$ we see that $\lim_{s \to 0} s\psi(s) = -1$ and hence we deduce that

(4.6) $\quad F'(0,v) = 1$

With $n=1$ in (4.3) we have

$$-\gamma_0(u) = -\frac{1}{2u} + \log u - \int_0^\infty e^{-uv}\left[\frac{1}{e^v-1} - \frac{1}{v} + \frac{1}{2}\right]dv$$

or equivalently, since $\psi(u) = -\gamma_0(u)$, we obtain the well-known integral

(4.7) $$\psi(u) = -\frac{1}{2u} + \log u - \int_0^\infty e^{-uv}\left[\frac{1}{e^v-1} - \frac{1}{v} + \frac{1}{2}\right]dv$$

This integral, which appears in [38, p.16], may be easily verified by differentiating Binet's first formula for $\log \Gamma(u)$ [39, p.249], which we saw above in (3.13).

We now consider the second derivative

$$F''(s,v) = -F(s,v)\psi'(s) + F(s,v)[\log v - \psi(s)]^2$$

$$= -F(s,v)\psi'(s) + F(s,v)\left[\log^2 v - 2\psi(s)\log v + \psi^2(s)\right]$$

and we have the limit

$$F''(0,v) = \lim_{s\to 0} F(s,v)\left[\psi^2(s) - \psi'(s)\right] - 2\log v \lim_{s\to 0} F(s,v)\psi(s)$$

$$= \lim_{s\to 0} F(s,v)\left[\left\{\psi^2(1+s) - 2\frac{\psi(1+s)}{s} + \frac{1}{s^2}\right\} - \left\{\psi'(1+s) + \frac{1}{s^2}\right\}\right] + 2\log v$$

$$= \lim_{s\to 0} F(s,v)\left[\psi^2(1+s) - 2\frac{\psi(1+s)}{s} - \psi'(1+s)\right] + 2\log v$$

Since $F(0,v) = 0$ this gives us

$$= -2\lim_{s\to 0} F(s,v)\left[\frac{\psi(1+s)}{s}\right] + 2\log v$$

and we therefore obtain

(4.8) $$F''(0,v) = 2(\gamma + \log v)$$

Referring back to (4.3) we have with $n=2$

$$\gamma_1(u) = \frac{1}{2u}\log u - \frac{1}{2}\log^2 u - \frac{1}{2}\int_0^\infty F^{(2)}(0,v)e^{-uv}\left[\frac{1}{e^v-1} - \frac{1}{v} + \frac{1}{2}\right]dv$$

and substituting (4.8) this becomes

$$\gamma_1(u) = \frac{1}{2u}\log u - \frac{1}{2}\log^2 u - \int_0^\infty (\gamma + \log v)e^{-uv}\left[\frac{1}{e^v-1} - \frac{1}{v} + \frac{1}{2}\right]dv$$

Using (4.7) we obtain

(4.9) $$\gamma_1(u) = \frac{1}{2u}\log u - \frac{1}{2}\log^2 u + \left[\psi(u) + \frac{1}{2u} - \log u\right]\gamma - \int_0^\infty e^{-uv}\left[\frac{1}{e^v-1} - \frac{1}{v} + \frac{1}{2}\right]\log v\, dv$$

and for $u = 1$ we have

(4.10) $$\gamma_1 = \left[-\gamma + \frac{1}{2}\right]\gamma - \int_0^\infty e^{-v}\left[\frac{1}{e^v-1} - \frac{1}{v} + \frac{1}{2}\right]\log v\, dv$$

Since $\Gamma^{(m)}(u) = \int_0^\infty v^{u-1} e^{-v} \log^m v\, dv$ we have $\Gamma^{(m)}(1) = \int_0^\infty e^{-v} \log^m v\, dv$ and, in particular, we see that $\Gamma'(1) = -\gamma = \int_0^\infty e^{-v} \log v\, dv$. Hence we may write (4.10) as

(4.11) $$\gamma_1 = -\gamma^2 + \gamma - \int_0^\infty e^{-v}\left[\frac{1}{e^v-1} - \frac{1}{v}\right]\log v\, dv$$

We note that

$$\frac{1}{1-e^{-v}} = \frac{e^v}{e^v-1} = 1 + \frac{1}{e^v-1}$$

with the result that

$$\int_0^\infty e^{-v}\left[\frac{1}{1-e^{-v}} - \frac{1}{v}\right]\log v\, dv = \int_0^\infty e^{-v}\left[\frac{1}{e^v-1} - \frac{1}{v}\right]\log v\, dv + \int_0^\infty e^{-v}\log v\, dv$$

$$= \int_0^\infty e^{-v}\left[\frac{1}{e^v-1} - \frac{1}{v}\right]\log v\, dv - \gamma$$

We then see that

(4.11.1) $$\gamma_1 = -\gamma^2 - \int_0^\infty e^{-v}\left[\frac{1}{1-e^{-v}} - \frac{1}{v}\right]\log v\, dv$$

as previously determined by Coppo [27] in 1999.

Having obtained integrals for $\gamma_0$ and $\gamma_1$, we now consider the general case for $\gamma_n$. We will approach this by reference to the (exponential) complete Bell polynomials, the salient features of which are summarised in the attached Appendix.

It is well known that [34]

$$(4.12) \qquad \frac{d^m}{dx^m} e^{f(x)} = e^{f(x)} Y_m\left(f^{(1)}(x), f^{(2)}(x), ..., f^{(m)}(x)\right)$$

where the (exponential) complete Bell polynomials $Y_n(x_1, ..., x_n)$ are defined by $Y_0 = 1$ and for $n \geq 1$

$$(4.13) \qquad Y_n(x_1, ..., x_n) = \sum_{\pi(n)} \frac{n!}{k_1! \, k_2! \, ... \, k_n!} \left(\frac{x_1}{1!}\right)^{k_1} \left(\frac{x_2}{2!}\right)^{k_2} ... \left(\frac{x_n}{n!}\right)^{k_n}$$

where the sum is taken over all partitions $\pi(n)$ of $n$, i.e. over all sets of integers $k_j$ such that

$$k_1 + 2k_2 + 3k_3 + ... + nk_n = n$$

For example, with $n = 1$ we see that the only possibility is $k_1 = 1$ and $k_j = 0 \; \forall \; j \geq 2$ which results in

$$Y_1(x_1) = x_1$$

With $n = 2$, we see that the possible outcomes are ($k_1 = 2$ and $k_2 = 0$) and ($k_1 = 0$ and $k_2 = 1$) which results in

$$Y_2(x_1, x_2) = x_1^2 + x_2$$

Suppose that $h'(x) = h(x)g(x)$ and let $f(x) = \log h(x)$. We see that

$$f'(x) = \frac{h'(x)}{h(x)} = g(x)$$

and then using (4.12) above we have

$$(4.14) \qquad \frac{d^m}{dx^m} h(x) = \frac{d^m}{dx^m} e^{\log h(x)} = h(x) Y_m\left(g(x), g^{(1)}(x), ..., g^{(m-1)}(x)\right)$$

As a variation of (4.14) above, suppose that $j'(x) = j(x)[g(x)+\alpha]$ where $\alpha$ is independent of $x$ and let $f(x) = \log j(x)$. We see that

$$f'(x) = \frac{j'(x)}{j(x)} = g(x) + \alpha$$

and $f^{(k+1)}(x) = g^{(k+1)}(x)$ for $k \geq 1$ and therefore we obtain

(4.15) $$\frac{d^m}{dx^m} j(x) = \frac{d^m}{dx^m} e^{\log j(x)} = j(x) Y_m\left(g(x)+\alpha, g^{(1)}(x),..., g^{(m-1)}(x)\right)$$

We note from (4.8) above that $F'(s,v) = F(s,v)[\log v - \psi(s)]$ and therefore we obtain

(4.16) $$F^{(n)}(s,v) = F(s,v) Y_n\left(\log v - \psi(s), -\psi^{(1)}(s),..., -\psi^{(n-1)}(s)\right)$$

We saw in (4.3) above that

(4.17) $$\gamma_{n-1}(u) = \frac{1}{2u}\log^{n-1} u - \frac{1}{n}\log^n u - \frac{(-1)^n}{n}\int_0^\infty F^{(n)}(0,v)e^{-uv}\left[\frac{1}{e^v-1} - \frac{1}{v} + \frac{1}{2}\right]dv$$

and hence we deduce that

(4.18)
$$(-1)^{n+1} n\gamma_{n-1} = \int_0^\infty F(s,v) Y_n\left(\log v - \psi(s), -\psi^{(1)}(s),..., -\psi^{(n-1)}(s)\right) e^{-v}\left[\frac{1}{e^v-1} - \frac{1}{v} + \frac{1}{2}\right]dv\bigg|_{s=0}$$

Prima facie, it is quite remarkable that this limit actually exists considering that the digamma function and all of its derivatives diverge at $s = 0$.

It is shown in Appendix A that

$$Y_n(x_1+\alpha, x_2,..., x_n) = \sum_{k=0}^n \binom{n}{k} \alpha^{n-k} Y_k(x_1,..., x_k)$$

and we then determine that

(4.19) $$\frac{d^m}{dx^m} j(x) = j(x) \sum_{k=0}^m \binom{m}{k} \alpha^{m-k} Y_k\left(g(x), g^{(1)}(x),..., g^{(k-1)}(x)\right)$$

Now, referring back to (4.5), we see that with $g(s) = \log v - \psi(s)$

(4.19.1) $$F^{(n)}(s,v) = F(s,v)\sum_{k=0}^{n}\binom{n}{k}\log^{n-k} v \cdot Y_k\left(-\psi(s),-\psi^{(1)}(s),...,-\psi^{(k-1)}(s)\right)$$

and thus we have

(4.20)
$$\gamma_{n-1}(u) = \frac{1}{2u}\log^{n-1} u - \frac{1}{n}\log^n u - \frac{(-1)^n}{n}\sum_{k=0}^{n}\binom{n}{k}\mathbf{Y}_k(s)\int_0^{\infty} F(s,v)\log^{n-k} v \cdot e^{-uv}\left[\frac{1}{e^v-1}-\frac{1}{v}+\frac{1}{2}\right]dv\bigg|_{s=0}$$

where, for convenience, we denote $\mathbf{Y}_k(s)$ as

$$\mathbf{Y}_k(s) = Y_k\left(-\psi(s),-\psi^{(1)}(s),...,-\psi^{(k-1)}(s)\right)$$

and we note that $\mathbf{Y}_k(s)$ is independent of the integration variable $v$.

$\square$

We now adopt a slightly different approach so as to eliminate all of the apparently troublesome factors $\psi^{(j)}(s)$ in the limit as $s \to 0$. To this end we write (4.2) in the equivalent form

$$F(s,v) = \frac{v^s}{\Gamma(s)} = \frac{sv^s}{\Gamma(1+s)}$$

and first of all we employ the Leibniz differentiation formula to obtain

$$\frac{\partial^n}{\partial s^n}F(s,v) = \sum_{k=0}^{n}\binom{n}{k}\frac{\partial^{n-k}}{\partial s^{n-k}}[sv^s]\frac{d^k}{ds^k}\frac{1}{\Gamma(1+s)}$$

We see that

$$\frac{d}{ds}\frac{1}{\Gamma(1+s)} = -\frac{1}{\Gamma(1+s)}\psi(1+s)$$

and applying (4.14) we determine that

(4.21) $$\frac{d^k}{ds^k}\frac{1}{\Gamma(1+s)} = \frac{1}{\Gamma(1+s)}Y_k\left(-\psi(1+s),-\psi^{(1)}(1+s),...,-\psi^{(k-1)}(1+s)\right)$$

We also have

$$\frac{\partial^j}{\partial s^j}[sv^s] = sv^s\log^j v + jv^s\log^{j-1} v$$

so that

$$\frac{\partial^n}{\partial s^n} F(s,v) = \frac{v^s}{\Gamma(1+s)} \sum_{k=0}^{n} \binom{n}{k} \left( s \log^{n-k} v + (n-k) \log^{n-k-1} v \right) \mathbf{Y}_k(1+s)$$

where, as before, we denote

$$\mathbf{Y}_k(1+s) = Y_k\left(-\psi(1+s), -\psi^{(1)}(1+s), \ldots, -\psi^{(k-1)}(1+s)\right)$$

When $s = 0$ this becomes

$$\frac{\partial^n}{\partial s^n} F(0,v) = \sum_{k=0}^{n} \binom{n}{k} (n-k) \log^{n-k-1} v\, Y_k\left(-\psi(1), -\psi^{(1)}(1), \ldots, -\psi^{(k-1)}(1)\right)$$

$$= \sum_{k=0}^{n-1} \binom{n}{k} (n-k) \log^{n-k-1} v\, Y_k\left(-\psi(1), -\psi^{(1)}(1), \ldots, -\psi^{(k-1)}(1)\right)$$

Using the elementary binomial identity [29, p.157]

$$(n-k)\binom{n}{k} = n \binom{n-1}{k}$$

this becomes

$$= n \sum_{k=0}^{n-1} \binom{n-1}{k} \log^{n-k-1} v\, Y_k\left(-\psi(1), -\psi^{(1)}(1), \ldots, -\psi^{(k-1)}(1)\right)$$

We therefore conclude that

(4.22) $\quad \gamma_{n-1}(u) = \frac{1}{2u} \log^{n-1} u - \frac{1}{n} \log^n u$

$$+ (-1)^{n-1} \sum_{k=0}^{n-1} \binom{n-1}{k} Y_k\left(-\psi(1), -\psi^{(1)}(1), \ldots, -\psi^{(k-1)}(1)\right) \int_0^\infty \log^{n-k-1} v \cdot e^{-uv} \left[\frac{1}{e^v - 1} - \frac{1}{v} + \frac{1}{2}\right] dv$$

and

(4.23)
$$\gamma_{n-1} = (-1)^{n-1} \sum_{k=0}^{n-1} \binom{n-1}{k} Y_k\left(-\psi(1), -\psi^{(1)}(1), \ldots, -\psi^{(k-1)}(1)\right) \int_0^\infty \log^{n-k-1} v \cdot e^{-v} \left[\frac{1}{e^v - 1} - \frac{1}{v} + \frac{1}{2}\right] dv$$

We note that

$$\int_0^\infty \log^{n-k-1} v \cdot e^{-v}\, dv = \Gamma^{(n-k-1)}(1)$$

and from (A.6) we have

$$\Gamma^{(k)}(1) = Y_k\left(\psi(1), \psi^{(1)}(1), \ldots, \psi^{(k-1)}(1)\right)$$

We will see in (A.9) that for $n \geq 1$ ($n = 0$ implies a value of 1)

(4.24) $$\sum_{j=0}^{n}\binom{n}{j} Y_j(x_1, \ldots, x_j) Y_{n-j}(-x_1, \ldots, -x_{n-j}) = \delta_{n,0}$$

and hence we may eliminate the factor of ½ in the integrand of (4.22) to obtain

(4.25) $$\gamma_{n-1}(u) = \frac{1}{2u}\log^{n-1} u - \frac{1}{n}\log^n u$$

$$+(-1)^{n-1}\sum_{k=0}^{n-1}\binom{n-1}{k} Y_k\left(-\psi(1), -\psi^{(1)}(1), \ldots, -\psi^{(k-1)}(1)\right)\int_0^\infty \log^{n-k-1} v \cdot e^{-uv}\left[\frac{1}{e^v - 1} - \frac{1}{v}\right] dv$$

(4.26)
$$\gamma_{n-1} = (-1)^{n-1}\sum_{k=0}^{n-1}\binom{n-1}{k} Y_k\left(-\psi(1), -\psi^{(1)}(1), \ldots, -\psi^{(k-1)}(1)\right)\int_0^\infty \log^{n-k-1} v \cdot e^{-v}\left[\frac{1}{e^v - 1} - \frac{1}{v}\right] dv$$

or equivalently for $n \geq 1$

(4.27) $$\gamma_n = (-1)^n \sum_{k=0}^{n}\binom{n}{k} Y_k\left(-\psi(1), -\psi^{(1)}(1), \ldots, -\psi^{(k-1)}(1)\right)\int_0^\infty \log^{n-k} v \cdot e^{-v}\left[\frac{1}{e^v - 1} - \frac{1}{v}\right] dv$$

We could also represent this in terms of the partial exponential Bell polynomials but this only seems to add an extra layer of complexity.

With the substitution $t = e^{-v}$ in (4.27) we obtain

(4.28)
$$\gamma_n = (-1)^n \sum_{k=0}^{n}\binom{n}{k} Y_k\left(-\psi(1), -\psi^{(1)}(1), \ldots, -\psi^{(k-1)}(1)\right)\int_0^1 \log^{n-k}[\log(1/t)]\left[\frac{t}{1-t} + \frac{1}{\log t}\right] dt$$

$$= (-1)^n \sum_{k=0}^{n}\binom{n}{k} Y_k\left(-\psi(1), -\psi^{(1)}(1), \ldots, -\psi^{(k-1)}(1)\right)\int_0^1 \log^{n-k}[\log(1/t)]\left[-1 + \frac{1}{1-t} + \frac{1}{\log t}\right] dt$$

Since $\int_0^1 \log^{n-k}(-\log(1/t))\,dt = \Gamma^{(n-k)}(1)$ we see, in the same manner as before, that for $n \geq 1$ the terms involving $-1$ in the integrand cancel out and we conclude that

$$(4.29) \quad \gamma_n = (-1)^n \sum_{k=0}^{n} \binom{n}{k} Y_k\left(-\psi(1), -\psi^{(1)}(1), \ldots, -\psi^{(k-1)}(1)\right) \int_0^1 \log^{n-k}[\log(1/t)] \left[\frac{1}{\log t} + \frac{1}{1-t}\right] dt$$

$\square$

In his 2001 dissertation, Brede [11] showed that there exists a polynomial $p_n(z)$ of degree $n$ such that

$$(4.30) \quad \gamma_n = \int_0^1 p_n[-\log\log(1/t)] \left[\frac{1}{\log t} + \frac{1}{1-t}\right] dt$$

Brede [11] stated, for example, that

$$p_0(z) = 1$$

$$p_1(z) = z - \gamma$$

$$p_2(z) = z^2 - 2\gamma z + \gamma^2 - \varsigma(2)$$

but he did not specify the precise form of the general polynomial $p_n(z)$.

Using these polynomials for $n = 1, 2, 3$ we have

$$\gamma = \int_0^1 \left[\frac{1}{\log t} + \frac{1}{1-t}\right] dt$$

$$\gamma_1 = -\int_0^1 \log[\log(1/t)] \left[\frac{1}{\log t} + \frac{1}{1-t}\right] dt - \gamma \int_0^1 \left[\frac{1}{\log t} + \frac{1}{1-t}\right] dt$$

so that

$$\int_0^1 \log[\log(1/t)] \left[\frac{1}{\log t} + \frac{1}{1-t}\right] dt = -\gamma_1 - \gamma^2$$

$$\gamma_2 = \int_0^1 \log^2[\log(1/t)] \left[\frac{1}{\log t} + \frac{1}{1-t}\right] dt + 2\gamma \int_0^1 \log[\log(1/t)] \left[\frac{1}{\log t} + \frac{1}{1-t}\right] dt$$

$$+[\gamma^2 - \varsigma(2)]\int_0^1 \left[\frac{1}{\log t} + \frac{1}{1-t}\right] dt$$

so that

$$\int_0^1 \log^2[\log(1/t)]\left[\frac{1}{\log t} + \frac{1}{1-t}\right] dt = 2\gamma\left[\gamma^2 + \gamma_1\right] - [\gamma^2 - \varsigma(2)]\gamma - \gamma_2$$

$\square$

Reindexing to $k = n - m$ we may write (4.29) as

$$\gamma_n = (-1)^n \sum_{m=n}^{n} \binom{n}{n-m} Y_{n-m}\left(-\psi(1), -\psi^{(1)}(1), \ldots, -\psi^{(n-m-1)}(1)\right) \int_0^1 \log^m[\log(1/t)]\left[\frac{1}{\log t} + \frac{1}{1-t}\right] dt$$

and since $\binom{n}{n-m} = \binom{n}{m}$ this becomes after reversing the order of summation

(4.31)
$$\gamma_n = \sum_{k=0}^{n} \binom{n}{k}(-1)^{n-k} Y_{n-k}\left(-\psi(1), -\psi^{(1)}(1), \ldots, -\psi^{(n-k-1)}(1)\right) \int_0^1 \left(-\log[\log(1/t)]\right)^k \left[\frac{1}{\log t} + \frac{1}{1-t}\right] dt$$

Comparing this with Brede's representation (4.30) we are therefore able to specify the precise form of Brede's polynomial

(4.32) $$p_n(z) = \sum_{k=0}^{n} \binom{n}{k}(-1)^{n-k} Y_{n-k}\left(-\psi(1), -\psi^{(1)}(1), \ldots, -\psi^{(n-k-1)}(1)\right) z^k$$

or equivalently

(4.33) $$p_n(z) = \sum_{k=0}^{n} \binom{n}{k}(-1)^k Y_k\left(-\psi(1), -\psi^{(1)}(1), \ldots, -\psi^{(k-1)}(1)\right) z^{n-k}$$

With regard to the above, we may note that

$$\frac{d^k}{ds^k} \frac{1}{\Gamma(s)} = Y_k\left(-\psi(s), -\psi^{(1)}(s), \ldots, -\psi^{(k-1)}(s)\right)$$

Differentiating (4.33) gives us

$$p_n'(z) = \sum_{k=0}^{n} (n-k)\binom{n}{k}(-1)^k Y_k\left(-\psi(1), -\psi^{(1)}(1), \ldots, -\psi^{(k-1)}(1)\right) z^{n-k-1}$$

$$= \sum_{k=0}^{n-1}(n-k)\binom{n}{k}(-1)^k Y_k\left(-\psi(1), -\psi^{(1)}(1), \ldots, -\psi^{(k-1)}(1)\right) z^{n-k-1}$$

$$= n\sum_{k=0}^{n-1}\binom{n-1}{k}(-1)^k Y_k\left(-\psi(1), -\psi^{(1)}(1), \ldots, -\psi^{(k-1)}(1)\right) z^{n-k-1}$$

and we therefore see that

(4.34) $\quad p_n'(z) = np_{n-1}(z)$

Since $p_n'(z) = np_{n-1}(z)$ we see that $p_n(z)$ is an Appell polynomial and therefore we have the relations

$$p_n(z) = \sum_{k=0}^{n}\binom{n}{k} p_{n-k}(0) z^k$$

$$p_n(x+y) = \sum_{k=0}^{n}\binom{n}{k} p_k(x) y^{n-k}$$

This concurs with Brede's result [11]

$$x^n = \int_0^\infty p_n(x - \log z) e^{-z} dz$$

We also note that

$$p_n'(z) = \sum_{k=0}^{n}\binom{n}{k}(-1)^{n-k} kY_{n-k}\left(-\psi(1), -\psi^{(1)}(1), \ldots, -\psi^{(n-k-1)}(1)\right) z^{k-1}$$

□

Having expended some energy getting to (4.22), it was somewhat disappointing to subsequently discover that this result could have been derived in a more succinct manner using (3.8). Subtracting a factor of $\dfrac{1}{s-1}$ from both sides of (3.8), we may write that equation in the following form

(4.50) $\quad \varsigma(s,u) - \dfrac{1}{s-1} = \dfrac{u^{-s}}{2} + \dfrac{u^{1-s}-1}{s-1} + \dfrac{1}{\Gamma(s)}\int_0^\infty e^{-uv} v^{s-1}\left[\dfrac{1}{e^v-1} - \dfrac{1}{v} + \dfrac{1}{2}\right] dv$

which we then differentiate and evaluate at $s=1$ this time to obtain

$$\frac{\partial^n}{\partial s^n}\left[\varsigma(s,u)-\frac{1}{s-1}\right]_{s=1} = \frac{(-1)^n}{2u}\log^n u + f^{(n)}(1) + \frac{\partial^n}{\partial s^n}\frac{1}{\Gamma(s)}\int_0^\infty e^{-uv}v^{s-1}\left[\frac{1}{e^v-1}-\frac{1}{v}+\frac{1}{2}\right]dv\bigg|_{s=1}$$

where, as a useful artifice, we have denoted $f(s)$ as

$$f(s) = \frac{u^{1-s}-1}{s-1}$$

We can represent $f(s)$ by the following integral

$$f(s) = \frac{u^{1-s}-1}{s-1} = -\int_1^u x^{-s}dx$$

so that

$$f^{(n)}(s) = -(-1)^n \int_1^u x^{-s} \log^n x\, dx$$

and thus we see that

$$f^{(n)}(1) = -(-1)^n \int_1^u \frac{\log^n x}{x} dx$$

$$= -(-1)^n \frac{\log^{n+1} u}{n+1}$$

Hence substituting (2.9)

$$\frac{\partial^n}{\partial s^n}\left[\varsigma(s,u)-\frac{1}{s-1}\right]_{s=1} = (-1)^n \gamma_n(u)$$

we obtain

$$(-1)^n \gamma_n(u) = \frac{(-1)^n}{2u}\log^n u - (-1)^n \frac{\log^{n+1} u}{n+1} + \frac{\partial^n}{\partial s^n}\frac{1}{\Gamma(s)}\int_0^\infty e^{-uv}v^{s-1}\left[\frac{1}{e^v-1}-\frac{1}{v}+\frac{1}{2}\right]dv\bigg|_{s=1}$$

Referring to the definition (4.2) of $F(s,v)$ we see that

$$\frac{\partial^n}{\partial s^n}\frac{1}{\Gamma(s)}\int_0^\infty e^{-uv}v^{s-1}\left[\frac{1}{e^v-1}-\frac{1}{v}+\frac{1}{2}\right]dv\bigg|_{s=1} = \frac{\partial^n}{\partial s^n}\int_0^\infty F(s,v)\frac{e^{-uv}}{v}\left[\frac{1}{e^v-1}-\frac{1}{v}+\frac{1}{2}\right]dv\bigg|_{s=1}$$

and using (4.19.1) this becomes

$$= \sum_{k=0}^{n}\binom{n}{k} Y_k\left(-\psi(s),-\psi^{(1)}(s),\ldots,-\psi^{(k-1)}(s)\right)\int_0^\infty e^{-uv}\log^{n-k} v\left[\frac{1}{e^v-1}-\frac{1}{v}+\frac{1}{2}\right]dv$$

Hence we obtain

(4.51) $\gamma_n(u) = \dfrac{1}{2u}\log^n u - \dfrac{1}{n+1}\log^{n+1} u$

$$+(-1)^n\sum_{k=0}^{n}\binom{n}{k} Y_k\left(-\psi(s),-\psi^{(1)}(s),\ldots,-\psi^{(k-1)}(s)\right)\int_0^\infty e^{-uv}\log^{n-k} v\left[\frac{1}{e^v-1}-\frac{1}{v}+\frac{1}{2}\right]dv$$

which corresponds with (4.22). Writing (4.50) as

(4.51.1) $\left[\varsigma(s,u)-\dfrac{1}{s-1}-\dfrac{u^{-s}}{2}-\dfrac{u^{1-s}-1}{s-1}\right]\Gamma(s) = \int_0^\infty e^{-uv} v^{s-1}\left[\dfrac{1}{e^v-1}-\dfrac{1}{v}+\dfrac{1}{2}\right]dv$

and, using the Leibniz formula to differentiate this, we obtain an inversion formula

(4.52)

$$\sum_{k=0}^{n}\binom{n}{k}(-1)^k\left[\gamma_k(u)-\dfrac{1}{2u}\log^k u+\dfrac{\log^{k+1} u}{k+1}\right]\Gamma^{(n-k)}(1) = \int_0^\infty e^{-uv}\log^n v\left[\dfrac{1}{e^v-1}-\dfrac{1}{v}+\dfrac{1}{2}\right]dv$$

and with $u = 1$ we have

(4.53) $\displaystyle\sum_{k=0}^{n}\binom{n}{k}(-1)^k \gamma_k \Gamma^{(n-k)}(1) - \dfrac{1}{2}\Gamma^{(n)}(1) = \int_0^\infty e^{-v}\log^n v\left[\dfrac{1}{e^v-1}-\dfrac{1}{v}+\dfrac{1}{2}\right]dv$

or equivalently

(4.54) $\displaystyle\sum_{k=0}^{n}\binom{n}{k}(-1)^k \gamma_k \Gamma^{(n-k)}(1) = \int_0^\infty e^{-v}\log^n v\left[\dfrac{1}{e^v-1}-\dfrac{1}{v}\right]dv$

Differentiating (4.51.1) with respect to $u$ gives us

$$\left[-s\varsigma(s+1,u)+\dfrac{su^{-s-1}}{2}+u^{-s}\right]\Gamma(s) = -\int_0^\infty e^{-uv} v^s\left[\dfrac{1}{e^v-1}-\dfrac{1}{v}+\dfrac{1}{2}\right]dv$$

which may be written as

$$\left[\varsigma(s+1)-\frac{1}{s}-\frac{1}{2}\right]\Gamma(1+s) = \int_0^\infty e^{-v}v^s\left[\frac{1}{e^v-1}-\frac{1}{v}+\frac{1}{2}\right]dv$$

and differentiating this will also result in (4.52)

(4.55) $\displaystyle\sum_{k=0}^n \binom{n}{k}(-1)^k\left[\gamma_k(u)-\frac{1}{2u}\log^k u+\frac{\log^{k+1} u}{k+1}\right]\Gamma^{(n-k)}(1) = \int_0^\infty e^{-uv}\log^n v\left[\frac{1}{e^v-1}-\frac{1}{v}+\frac{1}{2}\right]dv$

□

Referring back to the method originally employed by Stieltjes, in particular (3.22), we have

$$f_n(u) = \int_0^\infty \left[\frac{e^{-v}}{v}-\frac{e^{-(u+1)v}}{v}\right]\log^n v\, dv$$

$$f_n'(u) = F_n(u+1) = \sum_{k=0}^n (-1)^k \binom{n}{k}\Gamma^{(n-k)}(1)\frac{\log^k(u+1)}{u+1}$$

and since $f_n(0) = 0$ we have by integration

$$f_n(u) = \sum_{k=0}^n (-1)^k \binom{n}{k}\Gamma^{(n-k)}(1)\frac{\log^{k+1}(u+1)}{k+1}$$

so that

$$f_n(u-1) = \sum_{k=0}^n (-1)^k \binom{n}{k}\Gamma^{(n-k)}(1)\frac{\log^{k+1} u}{k+1}$$

$$\int_0^\infty \left[\frac{e^{-v}}{v}-\frac{e^{-uv}}{v}\right]\log^n v\, dv = \sum_{k=0}^n (-1)^k \binom{n}{k}\Gamma^{(n-k)}(1)\frac{\log^{k+1} u}{k+1}$$

Referring back to (4.52)

$$\sum_{k=0}^n \binom{n}{k}(-1)^k\left[\gamma_k(u)-\frac{1}{2u}\log^k u+\frac{\log^{k+1} u}{k+1}\right]\Gamma^{(n-k)}(1) = \int_0^\infty e^{-uv}\log^n v\left[\frac{1}{e^v-1}-\frac{1}{v}+\frac{1}{2}\right]dv$$

we see that

(4.56) $\displaystyle\sum_{k=0}^n \binom{n}{k}(-1)^k\left[\gamma_k(u)-\frac{1}{2u}\log^k u\right]\Gamma^{(n-k)}(1)$

$$= \int_0^\infty e^{-uv} \log^n v \left[ \frac{1}{e^v - 1} - \frac{1}{v} + \frac{1}{2} \right] dv - \int_0^\infty \left[ \frac{e^{-v}}{v} - \frac{e^{-uv}}{v} \right] \log^n v \, dv$$

$$= \int_0^\infty \left\{ e^{-uv} \left[ \frac{1}{e^v - 1} + \frac{1}{2} \right] - \frac{e^{-v}}{v} \right\} \log^n v \, dv$$

and with $u = 1$ we come back to (4.53) above.

$$\sum_{k=0}^n \binom{n}{k} (-1)^k \gamma_k \Gamma^{(n-k)}(1) - \frac{1}{2} \Gamma^{(n)}(1) = \int_0^\infty e^{-v} \log^n v \left[ \frac{1}{e^v - 1} - \frac{1}{v} + \frac{1}{2} \right] dv$$

With $n = 0, 1, 2$ we obtain

$$\gamma - \frac{1}{2} = \int_0^\infty e^{-v} \left[ \frac{1}{e^v - 1} - \frac{1}{v} + \frac{1}{2} \right] dv$$

$$-\gamma^2 - \gamma_1 + \frac{1}{2}\gamma = \int_0^\infty e^{-v} \log v \left[ \frac{1}{e^v - 1} - \frac{1}{v} + \frac{1}{2} \right] dv$$

$$\left(\gamma - \frac{1}{2}\right)[\gamma^2 + \varsigma(2)] + 2\gamma\gamma_1 + \gamma_2 = \int_0^\infty e^{-v} \log^2 v \left[ \frac{1}{e^v - 1} - \frac{1}{v} + \frac{1}{2} \right] dv$$

The approximate values of the first three Stieltjes constants are [5]

$$\gamma = 0.5772\cdots \qquad \gamma_1 = -0.0728\cdots \qquad \gamma_2 = -0.0096\cdots$$

and inserting these values into the above three equations numerically demonstrates that the corresponding three integrals have positive values.

□

We now wish to show that the integral $I_n$, which is defined below, is strictly positive for all positive integers and for $n = 0$

$$I_n = \int_0^\infty \log^n v \cdot e^{-v} \left[ \frac{1}{e^v - 1} - \frac{1}{v} + \frac{1}{2} \right] dv$$

We see that

$$I_n = \int_0^1 \log^n v \cdot e^{-v} \left[ \frac{1}{e^v - 1} - \frac{1}{v} + \frac{1}{2} \right] dv + \int_1^\infty \log^n v \cdot e^{-v} \left[ \frac{1}{e^v - 1} - \frac{1}{v} + \frac{1}{2} \right] dv$$

$$= J_n + K_n$$

With the substitution $t = e^{-v}$ we obtain for the first component $J_n$

$$\int_0^1 \log^n v \cdot e^{-v} \left[\frac{1}{e^v - 1} - \frac{1}{v} + \frac{1}{2}\right] dv = \int_0^{1/e} \log^n (\log(1/t)) \left[\frac{1}{1/t - 1} + \frac{1}{\log t} + \frac{1}{2}\right] dt$$

and with the substitution $t = 1/u$ this becomes

$$(4.56) \quad \int_0^{1/e} \log^n (\log(1/t)) \left[\frac{1}{1/t - 1} + \frac{1}{\log t} + \frac{1}{2}\right] dt = \int_e^\infty \log^n (\log u) \left[\frac{1}{u - 1} - \frac{1}{\log u} + \frac{1}{2}\right] \frac{du}{u^2}$$

It is clear that $\log^n (\log u) \geq 0$ for all $u \geq e$ and we now consider the other part of the integrand in (4.56). Let

$$\phi(u) = \frac{1}{u-1} - \frac{1}{\log u} + \frac{1}{2}$$

$$= \frac{(u+1)\log u - 2u + 2}{2(u-1)\log u}$$

and we note that $\phi(e) = \frac{1}{e-1} - \frac{1}{2} > 0$. The denominator of $\phi(u)$ is positive for all $u \geq e$ and we now consider the numerator

$$h(u) = (u+1)\log u - 2u + 2$$

We note that $h(e) = 3 - e > 0$ and we have the derivative

$$h'(u) = \log u + \frac{1}{u} - 1$$

and $h'(e) = e^{-1} > 0$. We see that

$$h''(u) = \frac{u-1}{u^2}$$

and therefore $h''(u) > 0$ for $u > 1$. This enables us to conclude that $h'(u)$ is monotonic increasing for $u > e$. Since $h'(e)$ is positive we then deduce that $h(u)$ is also monotonic increasing for $u > e$. Accordingly, $h(u) \geq 0$ for all $u \geq e$. Finally, we have $\phi(u) \geq 0$ for all $u \geq e$. Hence, since the integrand is positive, we conclude that $J_n$ is positive.

We now consider $K_n$ and we wish to determine the sign of $f(v)$ where

$$f(v) = \frac{1}{e^v - 1} - \frac{1}{v} + \frac{1}{2}$$

$$= \frac{2v + (v-2)(e^v - 1)}{2v(e^v - 1)}$$

The denominator of $f(v)$ is positive for all $v \geq 1$ and we now consider the numerator

$$g(v) = 2v + (v-2)(e^v - 1)$$

We note that $g(1) = 3 - e > 0$ and we have the derivative

$$g'(v) = (v-1)e^v + 1$$

so that $g'(v) > 0$ for $v \geq 1$. Accordingly, $g(v) \geq 0$ for all $v \geq 1$. Finally, we have $f(v) \geq 0$ for all $v \geq 1$. Hence, since the integrand is positive, we conclude that $K_n$ is positive.

It has therefore been demonstrated that $I_n = \int_0^\infty \log^n v \cdot e^{-v} \left[ \frac{1}{e^v - 1} - \frac{1}{v} + \frac{1}{2} \right] dv$ is positive for all positive integers and for $n = 0$.

Having regard to (3.7.1) we easily see that $\lim_{v \to 0} \left[ \frac{1}{e^v - 1} - \frac{1}{v} + \frac{1}{2} \right] = 0$; alternatively, this may be easily determined using L'Hôpital's rule.

It is possible that such representations for the Stieltjes constants may be useful in connection with the proof of the Riemann Hypothesis via the Li/Keiper constants (see for example [26]). This is because we note from (4.53) that

(4.57) $\sum_{k=0}^{n} \binom{n}{k} (-1)^k \gamma_k \Gamma^{(n-k)}(1) - \frac{1}{2} \Gamma^{(n)}(1) = \int_0^\infty e^{-v} \log^n v \left[ \frac{1}{e^v - 1} - \frac{1}{v} + \frac{1}{2} \right] dv > 0$

and reference to the Appendix then shows that if $n$ is even then $\Gamma^{(n)}(1)$ is positive and thus

(4.58) $\sum_{k=0}^{n} \binom{n}{k} (-1)^k \gamma_k \Gamma^{(n-k)}(1) > 0$

$\square$

Whilst not specifically employed in this paper, I accidentally came across the following lemma whilst I was trying to prove that $I_n$ was positive

Suppose that $f(t) \leq 0$ and $f'(t) \leq 0$ for all $t \in [1, \infty)$. We see that

$$I_n = \int_0^\infty f(t) \log^n t \, dt$$

$$= \int_0^1 f(t) \log^n t \, dt + \int_1^\infty f(t) \log^n t \, dt$$

and with the substitution $t = 1/y$ we have

$$\int_0^1 f(t) \log^n t \, dt = (-1)^n \int_1^\infty f(1/y) y^{-2} \log^n y \, dy$$

We therefore obtain

$$I_n = \int_1^\infty \left[ f(t) + (-1)^n t^{-2} f(1/t) \right] \log^n t \, dt$$

and it is clear that the integrand is negative in the interval $[1, \infty)$ in the case where $n$ is an even integer. We now consider the case where $n$ is an odd integer and we want to prove that the following integral is also negative

$$I_{2n+1} = \int_1^\infty \left[ f(t) - t^{-2} f(1/t) \right] \log^{2n+1} t \, dt$$

As a consequence of the Mean Value Theorem of calculus we have

$$f(t) - f(1/t) = \left( t - t^{-1} \right) f'(\alpha)$$

where $\alpha$ is such that $t > \alpha > t^{-1}$. We see that

$$f(t) - t^{-2} f(1/t) = f(t) - f(1/t) + f(1/t) - t^{-2} f(1/t)$$

$$= f(t) - f(1/t) + f(1/t)(1 - t^{-2})$$

$$= \left( t - \frac{1}{t} \right) f'(\alpha) + f(1/t)(1 - t^{-2})$$

$$= [t f'(\alpha) + f(1/t)](1 - t^{-2})$$

and this is clearly negative in the interval $(1, \infty)$. We therefore conclude that $I_{2n+1}$ is also negative.

We could also express $I_n$ as

$$I_n = \int_0^1 \left[ f(t) + (-1)^n t^{-2} f(1/t) \right] \log^n t \, dt$$

and hence $I_n$ will also be negative if $f(t) \leq 0$ and $f'(t) \leq 0$ for all $t \in (0,1]$.

**5. An application of the alternating Hurwitz zeta function**

The alternating Hurwitz zeta function $\varsigma_a(s, x)$ is defined by

$$\varsigma_a(s, x) = \sum_{n=0}^{\infty} \frac{(-1)^n}{(n+x)^s}$$

Upon a separation of terms according to the parity of $n$ we see that for $\text{Re}(s) > 1$

$$\varsigma_a(s, x) = \sum_{n=0}^{\infty} \frac{(-1)^n}{(n+x)^s} = \sum_{n=0}^{\infty} \frac{1}{(2n+x)^s} - \sum_{n=0}^{\infty} \frac{1}{(2n+1+x)^s}$$

$$= 2^{-s} \left[ \sum_{n=0}^{\infty} \frac{1}{(n+x/2)^s} - \sum_{n=0}^{\infty} \frac{1}{(n+(x+1)/2)^s} \right]$$

and we therefore see that $\varsigma_a(s, t)$ is related to the Hurwitz zeta function by the well-known formula [30]

(5.1) $$\varsigma_a(s, x) = 2^{-s} \left[ \varsigma\left(s, \frac{x}{2}\right) - \varsigma\left(s, \frac{1+x}{2}\right) \right]$$

Differentiation gives us

$$\frac{\partial}{\partial x} \varsigma_a(s, x) = -2^{-s-1} s \left[ \varsigma\left(s+1, \frac{x}{2}\right) - \varsigma\left(s+1, \frac{1+x}{2}\right) \right]$$

$$\frac{\partial^{n+1}}{\partial s^{n+1}} \frac{\partial}{\partial x} \varsigma_a(s, x) \bigg|_{s=0} = -\frac{1}{2}(n+1) \frac{\partial^n}{\partial s^n} [f(s, x) 2^{-s}] \bigg|_{s=0}$$

and using the Leibniz formula we obtain

$$= -\frac{1}{2}(n+1)\sum_{k=0}^{n}\binom{n}{k}(-1)^{n-k}2^{-s}\log^{n-k}2\cdot f^{(k)}(s,x)\bigg|_{s=0}$$

where $f(s,x) = \varsigma\left(s+1,\frac{x}{2}\right) - \varsigma\left(s+1,\frac{1+x}{2}\right)$

We have

$$f^{(k)}(0,x) = \varsigma^{(k)}\left(1,\frac{x}{2}\right) - \varsigma^{(k)}\left(1,\frac{1+x}{2}\right)$$

$$= (-1)^k\left[\gamma_k\left(\frac{x}{2}\right) - \gamma_k\left(\frac{1+x}{2}\right)\right]$$

and therefore

(5.2) $\quad \dfrac{\partial^{n+1}}{\partial s^{n+1}}\dfrac{\partial}{\partial x}\varsigma_a(s,x)\bigg|_{s=0} = \dfrac{1}{2}(-1)^{n+1}(n+1)\sum_{k=0}^{n}\binom{n}{k}\log^{n-k}2\cdot\left[\gamma_k\left(\dfrac{x}{2}\right) - \gamma_k\left(\dfrac{1+x}{2}\right)\right]$

We now refer to the Hasse identity for the alternating Hurwitz zeta function (see equation (4.4.79) in [21])

(5.3) $\quad \varsigma_a(s,x) = \sum_{i=0}^{\infty}\dfrac{1}{2^{i+1}}\sum_{j=0}^{i}\binom{i}{j}\dfrac{(-1)^j}{(x+j)^s}$

Differentiation gives us

$$\frac{\partial}{\partial x}\varsigma_a(s,x) = -s\sum_{i=0}^{\infty}\frac{1}{2^{i+1}}\sum_{j=0}^{i}\binom{i}{j}\frac{(-1)^j}{(x+j)^{s+1}}$$

$$\frac{\partial^{n+1}}{\partial s^{n+1}}\frac{\partial}{\partial x}\varsigma_a(s,x)\bigg|_{s=0} = (-1)^{n+1}(n+1)\sum_{i=0}^{\infty}\frac{1}{2^{i+1}}\sum_{j=0}^{i}\binom{i}{j}\frac{(-1)^j\log^n(x+j)}{x+j}$$

and equating this with (5.2) gives us

(5.4) $\quad \dfrac{1}{2}\sum_{k=0}^{n}\binom{n}{k}\log^{n-k}2\cdot\left[\gamma_k\left(\dfrac{x}{2}\right) - \gamma_k\left(\dfrac{1+x}{2}\right)\right] = \sum_{i=0}^{\infty}\dfrac{1}{2^{i+1}}\sum_{j=0}^{i}\binom{i}{j}\dfrac{(-1)^j\log^n(x+j)}{x+j}$

We have

$$\varsigma_a^{(n)}(s,x) = (-1)^n \sum_{i=0}^{\infty} \frac{(-1)^i \log^n(x+i)}{(x+i)^s} = (-1)^n \sum_{i=0}^{\infty} \frac{1}{2^{i+1}} \sum_{j=0}^{i} \binom{i}{j} \frac{(-1)^j \log^n(x+j)}{(x+j)^s}$$

and hence we may obtain expressions for $\varsigma_a^{(n)}(1,x)$ in terms of the Stieltjes constants

(5.5) $$\varsigma_a^{(n)}(1,x) = \frac{1}{2} \sum_{k=0}^{n} \binom{n}{k} \log^{n-k} 2 \cdot \left[ \gamma_k\left(\frac{x}{2}\right) - \gamma_k\left(\frac{1+x}{2}\right) \right]$$

With $x=1$ and $n=0$ in (5.4) we obtain

$$\frac{1}{2}\left[ \gamma_0\left(\frac{1}{2}\right) - \gamma_0(1) \right] = \sum_{i=0}^{\infty} \frac{1}{2^{i+1}} \sum_{j=0}^{i} \binom{i}{j} \frac{(-1)^j}{1+j}$$

which gives us

$$\log 2 = \sum_{i=0}^{\infty} \frac{1}{2^{i+1}} \sum_{j=0}^{i} \binom{i}{j} \frac{(-1)^j}{1+j}$$

This is equivalent to $\lim_{s \to 1} \varsigma_a(s) = \log 2$.

With $x=1$ and $n=1$ in (5.4) we get

$$\frac{1}{2}\left[ \gamma_0\left(\frac{1}{2}\right) - \gamma_0(1) \right] \log 2 + \frac{1}{2}\left[ \gamma_1\left(\frac{1}{2}\right) - \gamma_1 \right] = \sum_{i=0}^{\infty} \frac{1}{2^{i+1}} \sum_{j=0}^{i} \binom{i}{j} \frac{(-1)^j \log(1+j)}{1+j}$$

or

$$\log^2 2 + \frac{1}{2}\left[ \gamma_1\left(\frac{1}{2}\right) - \gamma_1 \right] = \sum_{i=0}^{\infty} \frac{1}{2^{i+1}} \sum_{j=0}^{i} \binom{i}{j} \frac{(-1)^j \log(1+j)}{1+j}$$

We have the relationship [23]

(5.6) $$\sum_{r=1}^{q-1} \gamma_p\left(\frac{r}{q}\right) = -\gamma_p + q(-1)^p \frac{\log^{p+1} q}{p+1} + q \sum_{j=0}^{p} \binom{p}{j} (-1)^j \gamma_{p-j} \log^j q$$

and with $q=2$ this becomes

(5.7) $$\gamma_p\left(\frac{1}{2}\right) = -\gamma_p + 2(-1)^p \frac{\log^{p+1} 2}{p+1} + 2 \sum_{j=0}^{p} \binom{p}{j} (-1)^j \gamma_{p-j} \log^j 2$$

and in particular we have

(5.8) $$\gamma_1\left(\frac{1}{2}\right) = -\gamma_1 - \log^2 2 + 2\gamma_1 - 2\gamma \log 2$$

We then obtain

(5.9) $$\gamma = \frac{1}{2}\log 2 - \frac{1}{\log 2}\sum_{i=0}^{\infty}\frac{1}{2^{i+1}}\sum_{j=0}^{i}\binom{i}{j}\frac{(-1)^j \log(1+j)}{1+j}$$

This expression for Euler's constant was originally derived by Coffey [17] in 2006 (a different derivation is contained in equation (4.4.116g) in [20])).

Similarly we may also obtain with $n = 2$

$$\gamma_1 = -\frac{1}{12}\log^2 2 + \frac{1}{2}\sum_{i=0}^{\infty}\frac{1}{2^{i+1}}\sum_{j=0}^{i}\binom{i}{j}(-1)^j \frac{\log(1+j)}{1+j} - \frac{1}{2\log 2}\sum_{i=0}^{\infty}\frac{1}{2^{i+1}}\sum_{j=0}^{i}\binom{i}{j}(-1)^j \frac{\log^2(1+j)}{1+j}$$

which was also previously determined by Coffey [17] in a different manner.

Letting $x = 1$ and $x = 2$ in (5.4) results in

$$\frac{1}{2}\sum_{k=0}^{n}\binom{n}{k}\log^{n-k} 2 \cdot \left[\gamma_k\left(\frac{1}{2}\right) - \gamma_k\right] = \sum_{i=0}^{\infty}\frac{1}{2^{i+1}}\sum_{j=0}^{i}\binom{i}{j}\frac{(-1)^j \log^n(1+j)}{1+j}$$

$$\frac{1}{2}\sum_{k=0}^{n}\binom{n}{k}\log^{n-k} 2 \cdot \left[\gamma_k - \gamma_k\left(\frac{3}{2}\right)\right] = \sum_{i=0}^{\infty}\frac{1}{2^{i+1}}\sum_{j=0}^{i}\binom{i}{j}\frac{(-1)^j \log^n(2+j)}{2+j}$$

We note from [38, p.89] that for $m \in \mathbb{N}$

$$\varsigma(s, m+x) = \varsigma(s, x) - \sum_{j=0}^{m-1}\frac{1}{(j+x)^s}$$

and thus

$$\frac{\partial^k}{\partial s^k}[\varsigma(s, m+x) - \varsigma(s, x)] = (-1)^{k+1}\sum_{j=0}^{m-1}\frac{\log^k(j+x)}{(j+x)^s}$$

Therefore, referring to (1.3) we obtain

$$\gamma_k(m+x) - \gamma_k(x) = -\sum_{j=0}^{m-1}\frac{\log^k(j+x)}{j+x}$$

and in particular we have

$$\gamma_k(1+x) - \gamma_k(x) = -\frac{\log^k x}{x}$$

Therefore we have

$$\gamma_k\left(\frac{3}{2}\right) - \gamma_k\left(\frac{1}{2}\right) = (-1)^{k+1} 2 \log^k 2$$

**APPENDIX**

**Some aspects of the (exponential) complete Bell polynomials**

It is well known that [33]

(A.1) $$\frac{d^m}{dx^m} e^{f(x)} = e^{f(x)} Y_m\left(f^{(1)}(x), f^{(2)}(x), \ldots, f^{(m)}(x)\right)$$

where the (exponential) complete Bell polynomials $Y_n(x_1, \ldots, x_n)$ are defined by $Y_0 = 1$ and for $n \geq 1$

(A.2) $$Y_n(x_1, \ldots, x_n) = \sum_{\pi(n)} \frac{n!}{k_1! \, k_2! \ldots k_n!} \left(\frac{x_1}{1!}\right)^{k_1} \left(\frac{x_2}{2!}\right)^{k_2} \cdots \left(\frac{x_n}{n!}\right)^{k_n}$$

where the sum is taken over all partitions $\pi(n)$ of $n$, i.e. over all sets of integers $k_j$ such that
$$k_1 + 2k_2 + 3k_3 + \ldots + nk_n = n$$

For example, with $n = 1$ we see that the only possibility is $k_1 = 1$ and $k_j = 0 \; \forall j \geq 2$ which results in

$$Y_1(x_1) = x_1$$

With $n = 2$, we see that the possible outcomes are ($k_1 = 2$ and $k_2 = 0$) and ($k_1 = 0$ and $k_2 = 1$) which results in

$$Y_2(x_1, x_2) = x_1^2 + x_2$$

Suppose that $h'(x) = h(x)g(x)$ and let $f(x) = \log h(x)$. We see that

$$f'(x) = \frac{h'(x)}{h(x)} = g(x)$$

and then using (A.1) above we have

(A.3) $$\frac{d^m}{dx^m} h(x) = \frac{d^m}{dx^m} e^{\log h(x)} = h(x) Y_m\left(g(x), g^{(1)}(x), \ldots, g^{(m-1)}(x)\right)$$

As an example, letting $h(x) = \Gamma(x)$ in (A.3) we obtain

(A.4) $$\frac{d^m}{dx^m} e^{\log \Gamma(x)} = \Gamma^{(m)}(x) = \Gamma(x) Y_m\left(\psi(x), \psi^{(1)}(x), \ldots, \psi^{(m-1)}(x)\right)$$

$$= \int_0^\infty t^{x-1} e^{-t} \log^m t \, dt$$

and since [26, p.22]

(A.5) $$\psi^{(p)}(x) = (-1)^{p+1} p! \varsigma(p+1, x)$$

we may express $\Gamma^{(m)}(x)$ in terms of $\psi(x)$ and the Hurwitz zeta functions. In particular, Kölbig [33] noted that

(A.6) $$\Gamma^{(m)}(1) = Y_m(-\gamma, x_1, \ldots, x_{m-1})$$

where $x_p = (-1)^{p+1} p! \varsigma(p+1)$. Values of $\Gamma^{(m)}(1)$ for $m \leq 10$ are reported in [38, p.265] and the first three are

$$\Gamma^{(1)}(1) = -\gamma$$

$$\Gamma^{(2)}(1) = \varsigma(2) + \gamma^2$$

$$\Gamma^{(3)}(1) = -[2\varsigma(3) + 3\gamma\varsigma(2) + \gamma^3]$$

As shown in [22], we note that $\Gamma^{(n)}(x)$ has the same sign as $(-1)^n$ for all $x \in (0, \alpha]$ where $\alpha$ is the unique positive root of $\psi(x) = 0$ (Gauss determined that $\alpha \approx 1.4616321\ldots$). This was also reported as an exercise in Apostol's book [3, p.303] for the particular case of $\Gamma^{(n)}(1)$.

We have from (A.5) and (4.21)

(A.6.1) $$\frac{d^k}{dx^k}\frac{1}{\Gamma(1+x)} = \frac{1}{\Gamma(1+x)}Y_k\left(-\psi(1+x), -1!\varsigma(2,1+x), \ldots, -(-1)^k(k-1)!\varsigma(k,1+x)\right)$$

so that with $x = 0$

(A.6.2) $$\left.\frac{d^k}{dx^k}\frac{1}{\Gamma(1+x)}\right|_{x=0} = Y_k\left(\gamma, -1!\varsigma(2), \ldots, -(-1)^k(k-1)!\varsigma(k)\right)$$

The complete Bell polynomials have integer coefficients and the first six are set out below [19, p.307]

(A.7) $$Y_1(x_1) = x_1$$

$$Y_2(x_1, x_2) = x_1^2 + x_2$$

$$Y_3(x_1, x_2, x_3) = x_1^3 + 3x_1x_2 + x_3$$

$$Y_4(x_1, x_2, x_3, x_4) = x_1^4 + 6x_1^2 x_2 + 4x_1 x_3 + 3x_2^2 + x_4$$

$$Y_5(x_1, x_2, x_3, x_4, x_5) = x_1^5 + 10x_1^3 x_2 + 10x_1^2 x_3 + 15x_1 x_2^2 + 5x_1 x_4 + 10x_2 x_3 + x_5$$

$$Y_6(x_1, x_2, x_3, x_4, x_5, x_6) = x_1^6 + 6x_1 x_5 + 15x_2 x_4 + 10x_3^2 + 15x_1^2 x_4 + 15x_2^3 + 60x_1 x_2 x_3$$

$$+ 20x_1^3 x_3 + 45x_1^2 x_2^2 + 15x_1^4 x_1 + x_6$$

The complete Bell polynomials are also given by the exponential generating function (Comtet [19, p.134])

(A.8) $$\exp\left(\sum_{j=1}^{\infty} x_j \frac{t^j}{j!}\right) = \sum_{n=0}^{\infty} Y_n(x_1, \ldots, x_n)\frac{t^n}{n!}$$

We see that

$$\left.\frac{d^n}{dt^n}\exp\left(\sum_{j=1}^{\infty} x_j \frac{t^j}{j!}\right)\right|_{t=0} = Y_n(x_1, \ldots, x_n)$$

We note that

$$\sum_{n=0}^{\infty} Y_n(ax_1, \ldots, ax_n)\frac{t^n}{n!} = \exp\left(\sum_{j=1}^{\infty} ax_j \frac{t^j}{j!}\right) = \exp a\left(\sum_{j=1}^{\infty} x_j \frac{t^j}{j!}\right) = \left[\exp\left(\sum_{j=1}^{\infty} x_j \frac{t^j}{j!}\right)\right]^a$$

and thus we have

$$\left[\sum_{n=0}^{\infty} Y_n(x_1,\ldots,x_n)\frac{t^n}{n!}\right]^a = \sum_{n=0}^{\infty} Y_n(ax_1,\ldots,ax_n)\frac{t^n}{n!}$$

We see with $a=-1$ that

$$1 = \sum_{n=0}^{\infty} Y_n(-x_1,\ldots,-x_n)\frac{t^n}{n!} \sum_{n=0}^{\infty} Y_n(x_1,\ldots,x_n)\frac{t^n}{n!}$$

and using the Cauchy product formula this becomes

$$= \sum_{n=0}^{\infty} \sum_{j=0}^{n} \binom{n}{j} Y_j(x_1,\ldots,x_j) Y_{n-j}(-x_1,\ldots,-x_{n-j})\frac{t^n}{n!}$$

Hence we deduce that for $n \geq 1$

(A.9) $$\sum_{j=0}^{n} \binom{n}{j} Y_j(x_1,\ldots,x_j) Y_{n-j}(-x_1,\ldots,-x_{n-j}) = 0$$

From the definition of the (exponential) complete Bell polynomials we have

$$Y_m(ax_1, a^2 x_2,\ldots, a^m x_m) = a^m Y_r(x_1,\ldots,x_m)$$

and thus with $a=-1$ we have

$$Y_m(-x_1, x_2,\ldots,(-1)^m x_m) = (-1)^m Y_m(x_1,\ldots,x_m)$$

We also note that

$$Y_m(x_1, -x_2,\ldots,(-1)^{m+1} x_m) = (-1)^m Y_m(-x_1,\ldots,-x_m)$$

but no discernible sign pattern emerges here.

We have the recurrence relation [35]

$$Y_{n+1}(x_1,\ldots,x_{n+1}) = \sum_{k=0}^{n} \binom{n}{k} Y_{n-k}(x_1,\ldots,x_{n-k}) x_{k+1} = \sum_{k=0}^{n} \binom{n}{k} Y_k(x_1,\ldots,x_k) x_{n-k+1}$$

$$= x_{n+1} + \sum_{k=1}^{n} \binom{n}{k} Y_k(x_1,\ldots,x_k) x_{n-k+1}$$

We note that

$$Y_n(x_1+y_1,...,x_n+y_n) = \exp\left(\sum_{j=1}^{\infty}(x_j+y_j)\frac{t^j}{j!}\right) = \exp\left(\sum_{j=1}^{\infty}x_j\frac{t^j}{j!}\right)\exp\left(\sum_{j=1}^{\infty}y_j\frac{t^j}{j!}\right)$$

$$= \sum_{n=0}^{\infty}Y_n(x_1,...,x_n)\frac{t^n}{n!} \cdot \sum_{n=0}^{\infty}Y_n(y_1,...,y_n)\frac{t^n}{n!}$$

and, as before, we apply the Cauchy series product formula to obtain

(A.10) $$Y_n(x_1+y_1,...,x_n+y_n) = \sum_{k=0}^{n}\binom{n}{k}Y_{n-k}(x_1,...,x_{n-k})Y_k(y_1,...,y_k)$$

and we note that

$$Y_k(\alpha,0,...,0) = \alpha^k$$

With $y_n = -x_n$ we obtain from (A.10)

(A.11) $$Y_n(0,...,0) = 0 = \sum_{k=0}^{n}\binom{n}{k}Y_{n-k}(x_1,...,x_{n-k})Y_k(-x_1,...,-x_k)$$

as in (A.9) above.

**REFERENCES**

[1] V.S.Adamchik, A Class of Logarithmic Integrals. Proceedings of the 1997 International Symposium on Symbolic and Algebraic Computation. ACM, Academic Press, 1-8, 2001.
http://www-2.cs.cmu.edu/~adamchik/articles/logs.htm

[2] V.S. Adamchik, On the Hurwitz function for rational arguments. Applied Mathematics and Computation, 187 (2007) 3-12.

[3] T.M. Apostol, Mathematical Analysis. Second Ed., Addison-Wesley Publishing Company, Menlo Park (California), London and Don Mills (Ontario), 1974.

[4] T.M. Apostol, Introduction to Analytic Number Theory. Springer-Verlag, New York, Heidelberg and Berlin, 1976.

[5] T.M. Apostol, Formulas for Higher Derivatives of the Riemann Zeta Function. Math. of Comp., 169, 223-232, 1985.

[6]  R.G. Bartle, The Elements of Real Analysis. 2$^{nd}$ Ed. John Wiley & Sons Inc, New York, 1976.

[7]  B.C. Berndt, The Gamma Function and the Hurwitz Zeta Function.
     Amer. Math. Monthly, 92, 126-130, 1985.

[8]  B.C. Berndt, Ramanujan's Notebooks. Part I, Springer-Verlag, 1985.

[9]  B.C. Berndt, Chapter eight of Ramanujan's Second Notebook.
     J. Reine Agnew. Math, Vol. 338, 1-55, 1983.
     http://www.digizeitschriften.de/no_cache/home/open-access/nach-zeitschriftentiteln/

[10] J. Bohman and C.-E. Fröberg, The Stieltjes Function-Definition and Properties.
     Math. of Computation, 51, 281-289, 1988.

[11] M. Brede, Eine reihenentwicklung der vervollständigten und ergänzten Riemannschen zetafunktion $\Xi(s) = s(s-1)\frac{\varsigma(s)\Gamma\left(\frac{s}{2}\right)}{\pi^{s/2}}$ $\left(s \in \mathbf{C} - \{0,1\}\right)$ und verwandtes.
     http://www.mathematik.uni-kassel.de/~koepf/Diplome/brede.pdf

[12] K. Chakraborty, S. Kanemitsu, and T. Kuzumaki, Finite expressions for higher derivatives of the Dirichlet $L-$function and the Deninger $R-$function.
     Hardy Ramanujan Journal, Vol.32 (2009) 38-53.
     http://59.90.235.217/hrj/kane4.pdf

[13] H. Chen, The Fresnel integrals revisited.
     Coll. Math. Journal, 40, 259-262, 2009.

[14] M.W. Coffey, Relations and positivity results for the derivatives of the Riemann $\xi$ function. J. Comput. Appl. Math., 166, 525-534 (2004)

[15] M.W. Coffey, New summation relations for the Stieltjes constants
     Proc. R. Soc. A, 462, 2563-2573, 2006.
     http://rspa.royalsocietypublishing.org/content/464/2091/711.full.pdf

[16] M.W. Coffey, New results on the Stieltjes constants: Asymptotic and exact evaluation. J. Math. Anal. Appl., 317 (2006) 603-612.
     math-ph/0506061 [abs, ps, pdf, other]

[17] M.W. Coffey, The Stieltjes constants, their relation to the $\eta_j$ coefficients, and representation of the Hurwitz zeta function. 2007.
     arXiv:math-ph/0612086 [ps, pdf, other]

[18] M.W. Coffey, On representations and differences of Stieltjes coefficients, and other relations. arXiv:0809.3277 [ps, pdf, other], 2008.

[19]  L. Comtet, Advanced Combinatorics, Reidel, Dordrecht, 1974.

[20]  D.F. Connon, Some series and integrals involving the Riemann zeta function, binomial coefficients and the harmonic numbers. Volume II(b), 2007. arXiv:0710.4024 [pdf]

[21]  D.F. Connon, Some series and integrals involving the Riemann zeta function, binomial coefficients and the harmonic numbers. Volume III, 2007. arXiv:0710.4025 [pdf]

[22]  D.F. Connon, Some applications of the Stieltjes constants. arXiv:0901.2083 [pdf], 2009.

[23]  D.F. Connon, New proofs of the duplication and multiplication formulae for the gamma and the Barnes double gamma functions. arXiv:0903.4539 [pdf], 2009.

[24]  D.F. Connon, The difference between two Stieltjes constants. arXiv:0906.0277 [pdf], 2009.

[25]  D.F. Connon, A recurrence relation for the Li/Keiper constants in terms of the Stieltjes constants. arXiv:0902.1691 [pdf], 2009.

[26]  D.F. Connon, Some possible approaches to the Riemann Hypothesis via the Li/Keiper constants. arXiv:1002.3484 [pdf], 2010.

[27]  M. A. Coppo, Nouvelles expressions des constantes de Stieltjes. Exposition. Math.17(4), 349–358, 1999.

[28]  C. Deninger, On the analogue of the formula of Chowla and Selberg for real quadratic fields. J. Reine Angew. Math., 351 (1984), 172–191. http://www.digizeitschriften.de/dms/toc/?PPN=PPN243919689_0351

[29]  R.L. Graham, D.E. Knuth and O. Patashnik, Concrete Mathematics. Second Ed. Addison-Wesley Publishing Company, Reading, Massachusetts, 1994.

[30]  E.R. Hansen and M.L. Patrick, Some Relations and Values for the Generalized Riemann Zeta Function. Math. Comput., Vol. 16, No. 79. (1962), pp. 265-274.

[31] H. Hasse, Ein Summierungsverfahren für Die Riemannsche $\varsigma$ - Reithe. Math. Z.32, 458-464, 1930. http://dz-srv1.sub.uni-goettingen.de/sub/digbib/loader?ht=VIEW&did=D23956&p=462

[32] C. Hermite, Correspondance d'Hermite et de Stieltjes. Gauthier-Villars, Paris, 1905. http://ebooks.library.cornell.edu/m/math/browse/author/a.php

[33] K.S. Kölbig, The complete Bell polynomials for certain arguments in terms of Stirling numbers of the first kind.
J. Comput. Appl. Math. 51 (1994) 113-116.
Also available electronically at:
A relation between the Bell polynomials at certain arguments and a Pochhammer symbol. CERN/Computing and Networks Division, CN/93/2, 1993.
http://doc.cern.ch/archive/electronic/other/preprints//CM-P/CM-P00065731.pdf

[34] K.S. Kölbig and W. Strampp, An integral by recurrence and the Bell polynomials. CERN/Computing and Networks Division, CN/93/7, 1993.
http://cdsweb.cern.ch/record/249027/

[35] J. Riordan, An introduction to combinatorial analysis. Wiley, 1958.

[36] R. Shail, A class of infinite sums and integrals.
Math. of Comput., Vol.70, No.234, 789-799, 2000.
A Class of Infinite Sums and Integrals

[37] R. Sitaramachandrarao, Maclaurin Coefficients of the Riemann Zeta Function.
*Abstracts Amer. Math. Soc.* **7**, 280, 1986.

[38] H.M. Srivastava and J. Choi, Series Associated with the Zeta and Related Functions. Kluwer Academic Publishers, Dordrecht, the Netherlands, 2001.

[39] E.T. Whittaker and G.N. Watson, A Course of Modern Analysis: An Introduction to the General Theory of Infinite Processes and of Analytic Functions; With an Account of the Principal Transcendental Functions. Fourth Ed., Cambridge University Press, Cambridge, London and New York, 1963.

[40] Z.X. Wang and D.R. Guo, Special Functions.
World Scientific Publishing Co Pte Ltd, Singapore, 1989.

[41] N.-Y. Zhang and K.S. Williams, Some results on the generalized Stieltjes constants.
Analysis 14, 147-162 (1994).
www.math.carleton.ca/~williams/papers/pdf/187.pdf

Donal F. Connon
Elmhurst
Dundle Road
Matfield
Kent TN12 7HD
dconnon@btopenworld.com